\documentclass[11pt]{amsart}

\usepackage{latexsym}
\usepackage{amssymb}
\usepackage{amsmath}

\newtheorem{theorem}{Theorem}[section]
\newtheorem{lemma}[theorem]{Lemma}
\newtheorem{proposition}[theorem]{Proposition}
\newtheorem{corollary}[theorem]{Corollary}
\newtheorem{definition}[theorem]{Definition}
\newtheorem{example}[theorem]{Example}
\newtheorem{remark}[theorem]{Remark}

\newcommand\Gr{\mathop{\rm Gr}}

\newcommand\supp{\mathop{\rm supp}}

\newcommand\tr{\mathop{\rm tr}}

\newcommand\whh{\mathop{\rm w^*h}}

\newcommand\nph{\varphi}

\newcommand\e{\epsilon}
\newcommand\T{\mathbb{T}}
\newcommand\R{\mathbb{R}}
\newcommand\op{\mathop{\rm op}}
\newcommand\loc{\mathop{\rm loc}}

\newcommand\nul{\mathop{\rm null}}
\newcommand\OL{\mathop{\rm OL}}

\newcommand{\cl}[1]{\mathcal{#1}}
\newcommand{\bb}[1]{\mathbb{#1}}

\begin{document}

\title{Closable multipliers}

\author{V.S. Shulman}
\address{Department of Mathematics, Vologda State Technical University, Vologda, Russia}
\email{shulman\_v@yahoo.com}

\author{I.G. Todorov}
\address{Department of Pure Mathematics, Queen's University Belfast, Belfast BT7 1NN, United Kingdom}
\email {i.todorov@qub.ac.uk}

\author{L. Turowska}
\address{Department of Mathematical Sciences,
Chalmers University of Technology and  the University of Gothenburg,
Gothenburg SE-412 96, Sweden}
\email{turowska@chalmers.se}


\date{January 25, 2010}

\begin{abstract}
Let $(X,\mu)$ and $(Y,\nu)$ be standard measure spaces.
A function $\nph\in L^\infty(X\times Y,\mu\times\nu)$ is called a (measurable)
Schur multiplier if the map $S_\nph$, defined on the space of Hilbert-Schmidt
operators from $L_2(X,\mu)$ to $L_2(Y,\nu)$ by multiplying their integral
kernels by $\nph$, is bound-ed in the operator norm.

The paper studies measurable functions $\nph$ for which
$S_\nph$ is closable in the norm topology or in the weak* topology.
We obtain a characterisation of w*-closable multipliers and relate
the question about
norm closability to the theory of operator synthesis.
We also study multipliers of two special types:
if $\nph$ is of Toeplitz type, that is, if $\nph(x,y)=f(x-y)$, $x,y\in G$,
where $G$ is a locally compact abelian group, then the closability of $\nph$
is related to the local inclusion of $f$
in the Fourier algebra $A(G)$ of $G$.  
If $\nph$ is a divided difference, that is, a function of the form
$(f(x)-f(y))/(x-y)$, then its closability is
related to the ``operator smoothness'' of the function $f$.
A number of examples of non-closable, norm closable
and w*-closable multipliers are presented.
\end{abstract}

\maketitle

\section{Introduction}

Let $(X,\mu)$ and $(Y,\nu)$ be standard measure spaces, $H_1
= L^2(X,\mu)$, $H_2 = L^2(Y,\nu)$, and let $\cl B(H_1,H_2)$ be the space
of all bounded linear operators from $H_1$ into $H_2$.
There is a method, due mainly to
Birman and Solomyak \cite{BS1,BS2, BS3, BS4}, to associate to certain
bounded measurable functions $\varphi$ on $X\times Y$, linear
transformations $S_{\varphi}$ on $\cl B(H_1,H_2)$; these
transformations are called {\it Schur multipliers} or, in the more
general setting of spectral measures in the place of $\mu$ and $\nu$, {\it double
operator integrals}. Namely, one first defines the map $S_{\varphi}$
on the space of all Hilbert-Schmidt operators by multiplying their
integral kernels by $\varphi$; if $S_{\varphi}$ is bounded in the
operator norm, one extends it to the space $\cl K(H_1,H_2)$ of all
compact operators by continuity. The map $S_{\varphi}$ is defined on
$\cl B(H_1,H_2)$ by taking the second dual of the constructed map on
$\cl K(H_1,H_2)$.

The function $\varphi$ is thus called a Schur multiplier if the
multiplication map $S_{\varphi}$, initially defined on the space of
Hilbert-Schmidt operators, is bounded in the operator norm.
Equivalently, $\varphi$ is a Schur multiplier if
$\varphi(x,y)k(x,y)$ is the integral kernel of a nuclear operator
provided $k(x,y)$ is such. The set of all Schur multipliers will be denoted by $\frak{S}(X,Y)$.

Many years ago, Professor Solomyak informed the first author that at
the early stages of the development of the theory of double operator integrals, there existed an
idea to define $S_{\varphi}$ for a more general class of functions
$\varphi$ (perhaps for all measurable ones) as the closure of a
densely defined linear operator in the operator norm, or in the weak
operator, topology. However, this approach was not pursued because no
information on the closability of the multiplication maps
$S_{\varphi}$ was obtained at that time.

The aim of this paper is to study the classes of functions $\varphi$
for which $S_{\varphi}$ is closable in the norm topology, or in the
weak* topology, of $\cl B(H_1,H_2)$. By weak* topology we mean the topology
$\sigma(\cl B(H_1,H_2),\cl C_1(H_2,H_1))$ induced on $\cl
B(H_1,H_2)$ by its duality with the space $\cl C_1(H_2,H_1)$
of nuclear operators from $H_2$ into $H_1$.
We denote these classes of functions by $\frak{S}_{{\rm cl}}(X,Y)$ and $\frak{S}_{w^*}(X,Y)$, respectively.

We obtain a satisfactory characterisation of the class $\frak{S}_{w^*}(X,Y)$.
In order to describe it,
let us denote by $\frak{S}_{\loc}(X,Y)$ the class of
functions $\varphi$ with the following property: for each
$\varepsilon>0$ there exist subsets $X_{\varepsilon}\subseteq X$ and
$Y_{\varepsilon}\subseteq Y$ of measure not exceeding $\varepsilon$, such that
the restriction of $\varphi$ to $(X\setminus X_{\varepsilon})\times
(Y\setminus Y_{\varepsilon})$ is a Schur multiplier.  We prove in
Theorem \ref{local} that the elements of  $\frak{S}_{\loc}(X,Y)$ (which we call local
Schur multipliers) can be characterised in terms similar to those of Peller's
characterisation of Schur multipliers \cite{peller}. Namely, they are precisely the functions of
the form $(a(x),b(y))$, where $a(\cdot)$ (resp. $b(\cdot)$) is
a measurable function from $X$ (resp. $Y$) into a
separable Hilbert space. Then we show in Theorem \ref{wstarmult}
that the w*-closable multipliers (that is, the elements of $\frak{S}_{w^*}(X,Y)$) are
precisely the functions of the form $t(x,y)/s(x,y)$ where $t$ and $s$ are
local Schur multipliers and $s(x,y) \neq 0$ for marginally almost all $(x,y)$.
In particular, $\frak{S}_{w^*}(X,Y)$ is an algebra of (equivalent classes of) functions.

For any measurable function $\varphi$ on $X\times Y$, there exists a maximal (in
a sense that we make precise in Section \ref{s_wstar}) countable family of rectangles on
each of which $\varphi$ is w*-closable; the complement of their union is
denoted by $\kappa_{\varphi}^{w^*}$. The \lq\lq size'' of
$\kappa_{\varphi}^{w^*}$ can be considered as a measure of the
extent to which $\varphi$ fails to be a w*-closable multiplier.

The information we obtain about $\frak{S}_{{\rm cl}}(X,Y)$ is less
precise. Roughly speaking, we show that, in order to verify whether
$\nph$ belongs to $\frak{S}_{{\rm cl}}(X,Y)$, one needs to check whether the set
$\kappa_{\varphi}^{w^*}$ supports a non-zero compact operator.
More precisely: if $\kappa_{\varphi}^{w^*}$ does not support a non-zero compact
operator then  $\nph\in \frak{S}_{{\rm cl}}(X,Y)$ while if there exists a
non-zero compact operator in the {\it smallest} masa-bimodule with
support $\kappa_{\varphi}^{w^*}$ then $\nph\notin
\frak{S}_{{\rm cl}}(X,Y)$. The difference between the smallest
masa-bimodule with support $\kappa_{\varphi}^{w^*}$ and the largest one
(the set of all operators supported on $\kappa_{\varphi}^{w^*}$) is
subtle; it is the subject of the theory of operator synthesis
\cite{arveson, sht}, an operator analogue of the theory of spectral
synthesis \cite{graham, katznelson}.  We prove that
$\frak{S}_{{\rm cl}}(X,Y)$ is an algebra and present various examples of
multipliers which are not norm-closable and of norm-closable multipliers which are not
w*-closable.

The product of two measure spaces possesses natural
\lq\lq pseudo-topological structures'', namely the $\omega$-topology and the $\tau$-topology,
which are related to the problem of closability of multipliers. A
set is called $\tau$-open (resp. $\omega$-open) if it is a countable union
of measurable rectangles and a null set (resp. a
set contained in $(X_0\times Y) \cup (X\times Y_0)$, where $X_0$ and
$Y_0$ are null sets). Denoting by $C_{\tau}(X\times Y)$
(resp. $C_{\omega}(X\times Y)$) the space of all $\tau$-continuous (resp. $\omega$-continuous)
complex valued functions on $X\times Y$, we prove that
$$\frak{S}_{{\rm cl}}(X,Y)\subseteq C_{\tau}(X\times Y) \text{   and   } \frak{S}_{w^*}(X,Y)\subseteq C_{\omega}(X\times Y)$$
(if one identifies functions equivalent with respect to the product measure).
Both inclusions are shown to be strict.

We present examples which show that in the chain
$$ \frak{S}(X,Y)\subseteq \frak{S}_{\loc}(X,Y)\subseteq \frak{S}_{w^*}(X,Y)\subseteq \frak{S}_{{\rm cl}}(X,Y)$$
the first and the third inclusions are strict. The question of whether the second inclusion is strict is left open.

The paper is organised as follows. In Section 2 we state some basic
definitions and results about subsets of, and functions on, product 
measure spaces, bimodules over maximal abelian selfadjoint algebras, Schur multipliers and
closable operators.
In Section 3 we examine local Schur multipliers, while Sections 4 and 5 are
devoted to the study of
w*-closable and norm-closable multipliers, respectively.

In Sections 6 and 7 we study multipliers of special types. Given
a complex function $f$ defined on a subinterval of the real
line, one may consider its {\it divided difference}, in other words, the
function $\check f$ on two variables given by $\check f (x,y)
= (f(x)-f(y))/(x-y)$. The corresponding class of multipliers plays
an important role in Perturbation Theory and Spectral Theory (see,
for example, \cite{peller} and the references therein). We show that
such multipliers are always norm-closable; in Theorem \ref{kf} and
Corollary \ref{div-loc} we formulate necessary and sufficient
conditions for $\check f$ to be a local Schur multiplier.

{\it Toeplitz multipliers} are (Schur, local, norm-closable or w*-closable)
multipliers $\varphi$ of the form
$\varphi(x,y) = f(x-y)$, where $f$ is a complex function defined on
a locally compact abelian group $G$ (equipped with Haar measure $\mu$). Theorem \ref{toeplitz} asserts
that a function  $f(x-y)$ is a w*-closable multiplier if and only if it
is a local Schur multiplier, and that this occurs precisely when $f$
is equivalent to a function which  belongs locally to the Fourier
algebra of $G$.
The closability of Toeplitz multipliers is shown to be related to
some questions about sets of multiplicity in harmonic
analysis.
We describe also those functions of the form $f(x-y)$ which are
($(\mu\times \mu)$-equivalent to) $\omega$-continuous or $\tau$-continuous functions.
En route, we obtain an example of a continuous (hence $\omega$-continuous)
function on $G\times G$ which is not a norm-closable multiplier.

\medskip

{\bf Acknowledgements.} The authors are very grateful to M.
Roginskaya and T. W. K${\rm \ddot{o}}$rner for their friendly and
helpful advice.

The third author was supported by the Swedish Research Council. 

\section{ Preliminary results}

\subsection{Pseudo-topologies on the product of measure spaces}\label{sec2}

In what follows, we write $\frak{B}(Z) = \frak{B}(Z,\gamma)$ for the algebra,
with respect to the pointwise product, of all
measurable complex valued functions defined on a measure space
$(Z,\gamma)$.
Let $(X,\mu)$ and $(Y,\nu)$ be standard $\sigma$-finite
measure spaces, which will be fixed throughout the paper.
We will often write $\frak{B}(X\times Y)$ in the place of
$\frak{B}(X\times Y,\mu\times\nu)$.

A subset of $X\times Y$ is said to be a
{\it measurable rectangle} (or simply a {\it rectangle})
if it is of the form
$\alpha\times\beta$, where $\alpha\subseteq X$ and $\beta\subseteq Y$
are measurable subsets. A subset
$E\subseteq X\times Y$ is called {\it marginally null} if
$E\subseteq (X_0\times Y)\cup(X\times Y_0)$, where
$\mu(X_0) = \nu(Y_0) = 0$. We call two subsets $E,F\subseteq X\times Y$
{\it marginally equivalent} (and write $E\simeq F$) if the symmetric difference
of $E$ and $F$ is marginally null. We say that $E$ {\it marginally contains} $F$
(or $F$ {\it is marginally contained in} $E$) if $F\setminus E$ is marginally null;
 $E$ and $F$ are said to be marginally disjoint if $E\cap F$ is marginally null.

A subset $E$ of $X\times Y$ is called {\it $\omega$-open}
if it is marginally equivalent to the union of a
countable set of rectangles. The complements of
$\omega$-open sets are called {\it $\omega$-closed}.
It is clear that the class of all $\omega$-open
(resp. $\omega$-closed) sets is closed under countable unions (resp.
intersections) and finite intersections (resp. unions); in other words,
the $\omega$-open sets form a pseudo-topology.

The following lemma shows that in some cases
one can form a certain kind of a union of a given,
possibly uncountable, family of $\omega$-open subsets.

\begin{lemma}\label{union}
Let $\mathcal{E}$ be a family of $\omega$-open subsets of $X\times Y$. Let $\mathcal{E}_{\sigma}$ be the set of all
countable unions of elements of $\mathcal{E}$. Then there exists a (unique up to marginal
equivalence) set $E\in \mathcal{E}_{\sigma}$ which
marginally contains every set in $\mathcal{E}$.
\end{lemma}
\begin{proof} First assume that the measures $\mu$ and $\nu$ are finite. On the set $\Pi$ of all measurable
rectangles we introduce a metric $\rho$, setting, for $R_1 = X_1\times Y_1$, $R_2 = X_2\times Y_2$ in $\Pi$,
$\rho(R_1,R_2) = \mu(X_1\triangledown X_2) + \nu(Y_1\triangledown Y_2)$ (here $\triangledown$
denotes symmetric difference). Then $\Pi$ is a separable metric space
whence the set
$\mathcal{F}$ of all rectangles that are contained in elements of $\mathcal{E}$ is also separable. Let $\{R_n: n\ge
1\}$ be a dense sequence in $\mathcal{F}$ and $E = \cup_{n=1}^{\infty} R_n$. Then it is clear that $E$ marginally contains
all $R\in \mathcal{F}$ and therefore all $E\in \mathcal{E}$.

In the general case, let $X = \cup_{n=1}^{\infty} X_n$ and $Y
= \cup_{m=1}^{\infty} Y_m$ with $\mu(X_n) < \infty$ and $\nu(Y_m) < \infty$, $n,m\in \bb{N}$.
For each pair $n,m$,
let $\cl E_{n,m} = \{E\cap (X_n\times Y_m) : E\in \mathcal{E}\}$, let $E_{n,m}$ be the
$\omega$-union of $\cl E_{n,m}$, and set $E = \cup_{n,m=1}^{\infty} E_{n,m}$.

The uniqueness is obvious.
\end{proof}

The set whose existence is guaranteed by Lemma \ref{union}
will be called the \emph{$\omega$-union} of $\mathcal{E}$.

We will say that two functions $\varphi,\psi\in \frak{B}(X\times Y)$ are equivalent,
and write $\varphi\sim \psi$, if the set
$D = \{(x,y)\in X\times Y : \varphi(x,y)\neq \psi(x,y)\}$ is null with respect to the product measure. If $D$
is marginally null then we say that $\varphi$ and $\psi$ coincide marginally everywhere, or that they are marginally
equivalent, and write $\varphi \simeq \psi$.
By $L^{\infty}(X\times Y)$ we denote as usual the subalgebra of $\frak{B}(X\times Y)$
of all (equivalence classes of) essentially bounded functions.

Let $E\subseteq X\times Y$ be an $\omega$-open set.
A measurable function $\nph : E\rightarrow \bb{C}$ is called
{\it $\omega$-continuous} if the set $\nph^{-1}(G)$ is $\omega$-open for every open subset $G\subseteq\bb{C}$.
As in \cite[Corollary 3.2]{eks} one can see that
the set $C_{\omega}(E)$ of all $\omega$-continuous functions on $E$ is a
subalgebra of $\frak{B}(E)$.
For an arbitrary set $\mathcal{M}\subseteq C_{\omega}(X\times Y)\subseteq \frak{B}(X\times Y) $,
we let its {\it null set} $\nul(\mathcal{M})$ be the complement of
the $\omega$-union of the family ${\mathcal E} = \{h^{-1}(\mathbb{C}\setminus \{0\}): h\in \mathcal{M}\}$.

We will need the following simple result from \cite{ks} (see the remark after \cite[Proposition 8.1]{ks}).

\begin{lemma}\label{wopen}
Let $E\subseteq X\times Y$ be an $\omega$-open set and let
$f : E\rightarrow\bb{C}$ be an $\omega$-continuous function.
If $f\sim 0$ then $f\simeq 0$.
\end{lemma}

Thus if two $\omega$-continuous functions are equivalent then they are marginally equivalent, and if a function is equivalent to an $\omega$-continuous function
then the latter is defined uniquely up to marginal equivalence.

We will also need another pseudo-topology on $X\times Y$.
Two subsets $E_1$, $E_2$ of $X\times Y$ will be called $\mu\times \nu$-equivalent if
their symmetric difference is a  $\mu\times \nu$-null set.
We will say that a subset $E\subseteq X\times Y$ is {\it $\tau$-open}
if it is $\mu\times \nu$-equivalent to a countable union of rectangles.
It is clear that the pseudo-topology $\tau$ is stronger than $\omega$.

\subsection{Bimodules}

If $H_1$ and $H_2$ are Hilbert spaces,
we denote by $\cl B(H_1,H_2)$ the space of all bounded linear operators
from $H_1$ into $H_2$, and by $\cl K(H_1,$ $H_2)$
(resp. $\cl C_1(H_1,H_2)$, $\cl C_2(H_1,H_2)$) the
space of compact (resp. nuclear, Hilbert-Schmidt) operators in $\cl B(H_1,H_2)$.
Let $\|T\|_{\op}$ denote the operator norm of $T\in \cl B(H_1,H_2)$.
As usual, we write $\cl B(H) = \cl B(H,H)$.
The space $\cl C_1(H_2,H_1)$ (resp. $\cl B(H_1,H_2)$)
can be naturally identified with the Banach space dual of
$\cl K(H_1,H_2)$ (resp. $\cl C_1(H_2,H_1)$),
the duality being given by the map
$(T,S)\to \langle T,S\rangle \stackrel{def}{=}\tr(TS)$.
Here $\tr A$ denotes the trace of a nuclear operator $A$.
For a subset $\cl W\subseteq \cl C_1(H_2,H_1)$,
let $\cl W^{\perp} = \{T\in \cl B(H_1,H_2) :
\langle T,S\rangle = 0, \ \mbox{ for each} \ S\in \cl W\}$.

For the rest of the paper we let $H_1 = L^2(X,\mu)$ and $H_2 = L^2(Y,\nu)$.
The space $L^2(X\times Y)$ will be identified with
$\cl C_2(H_1,H_2)$ via the map sending an element $k\in L^2(X\times Y)$
to the integral operator $I_k$ given by $I_k \xi (y) = \int_{X} k(x,y)\xi(x) d\mu(x)$,
$\xi\in H_1$, $y\in Y$.
In a similar fashion, $\cl C_1(H_2,H_1)$ will be identified with
the space $\Gamma(X,Y)$ of all functions $F : X\times Y\to{\mathbb C}$
which admit a representation
$$F(x,y)=\sum_{i=1}^{\infty} f_i(x)g_i(y),$$
where $f_i\in L^2(X,\mu)$, $g_i\in L^2(Y,\nu)$, $i\in \bb{N}$,
$\sum_{i=1}^{\infty}\|f_i\|^2_{2} < \infty$ and
$\sum_{i=1}^{\infty}\|g_i\|^2_{2}<\infty$. Equivalently, $\Gamma(X,Y)$
can be defined as the projective tensor product
$L^2(X,\mu)\hat\otimes L^2(Y,\nu)$.
It was shown in \cite[Theorem 6.5]{eks} that $\Gamma(X,Y)$ consists of $\omega$-continuous functions.
For brevity, we often identify a function $h\in \Gamma(X,Y)$ with the
corresponding integral operator $I_h\in \cl C_1(H_2,H_1)$. It will
be convenient to write $\Gamma(X\times Y)$ in the place of
$\Gamma(X,Y)$ (this allows for example to write $\Gamma(\kappa)$, where $\kappa$ is a rectangle).

If $f\in L^{\infty}(X,\mu)$, let $M_f\in \cl B(H_1)$
denote the operator of multiplication by $f$.
We will often identify the collection $\{M_f : f\in L^{\infty}(X,\mu)\}$ of all such operators with
the function space $L^{\infty}(X,\mu)$.
If $\alpha\subseteq X$ is measurable, we write $P(\alpha) = M_{\chi_{\alpha}}$
for the multiplication by the characteristic function of the set $\alpha$.
Similar definitions and identifications are made for $L^{\infty}(Y,\nu)$.
A subspace
$\cl W\subseteq \cl B(H_1,H_2)$ will be called a {\it bimodule} if $M_{\psi}TM_{\nph} \in \cl  W$
for all $T\in \cl W$,
$\nph\in L^{\infty}(X,\mu)$ and $\psi\in L^{\infty}(Y,\nu)$.
One defines bimodules in $\cl B(H_2,H_1)$ in a similar fashion.

We say that an $\omega$-closed subset $\kappa\subseteq X\times Y$
{\it supports} an operator $T\in \cl B(H_1,H_2)$
(or that $T$ {\it is supported on} $\kappa$)
if $P(\beta)TP(\alpha) = 0$ whenever $\alpha\times\beta$ is
marginally disjoint from $\kappa$.
For any subset ${\mathcal M}\subseteq \cl B(H_1,H_2)$, there exists a smallest
(up to marginal equivalence) $\omega$-closed set
$\supp\cl M$ which supports every operator $T\in{\mathcal M}$ \cite{eks}.
By \cite{sht}, for any $\omega$-closed set $\kappa$ there exists a smallest (resp. largest)
w*-closed bimodule ${\mathfrak M}_{\min}(\kappa)$ (resp. ${\mathfrak M}_{\max}(\kappa)$) with support $\kappa$ in the sense that
if $\mathfrak M\subseteq \cl B(H_1,H_2)$ is a w*-closed bimodule with
$\text{supp }{\mathfrak M} = \kappa$ then
${\mathfrak M}_{\min}(\kappa)\subseteq {\mathfrak M}\subseteq {\mathfrak M}_{\max}(\kappa)$.
If ${\mathfrak M}_{\min}(\kappa) = {\mathfrak M}_{\max}(\kappa)$, the set $\kappa$ is called {\it synthetic}.
By \cite[Theorem~4.4]{sht},
${\mathfrak M}_{\min}(\kappa) = \{I_h : h \in \Gamma(X,Y), h\chi_{\kappa} \simeq 0\}^{\perp}$.


\begin{lemma}\label{sequ}
Let $\cl W\subseteq \cl C_1(H_2,H_1)$ be a bimodule and
 $\{f_n\}_{n=1}^{\infty}\subseteq \Gamma(X,$ $Y)$ be a sequence such that
$\{I_{f_n}\}_{n=1}^{\infty}$ is dense in $\cl W$. Then
$\nul(\cl W) = \supp(\cl W^{\perp})$ $=$ $\cap_{n=1}^{\infty}f_n^{-1}(0)$.

In particular, $\cl W$ is norm dense in $\cl C_1(H_2,H_1)$ if and only if
there exists a sequence $\{h_n\}_{n=1}^{\infty}\subseteq \Gamma(X,Y)$
such that $I_{h_n}\in \cl W$ for every $n\in \bb{N}$
and the set $\cap_{n=1}^{\infty} h_n^{-1}(0)$ is marginally null.
\end{lemma}
\begin{proof}
We start by showing the second statement.
Since the Hilbert spaces $H_1$ and $H_2$ are separable,
the space $\cl C_1(H_2,H_1)$ is separable and hence there exists a sequence $\{K_n\}_{n=1}^{\infty}$
dense in $\cl W$. Suppose that $K_n = I_{h_n}$, where $h_n\in \Gamma(X,Y)$ and let
$E = \cap_n h_n^{-1}(0)$. If $E$ is not marginally null
then,
by \cite[Corollary 4.1]{sht}, ${\mathfrak M}_{\min}(E)$ contains a
non-zero operator $T$. By \cite[Theorem 4.4]{sht}, $\langle T, K_n\rangle = 0$
for all $n\in \bb{N}$, and hence $\cl W$ is not dense.

Conversely, suppose that $\cl W$ is not dense,  let $h_n\in \Gamma$, $n\in \bb{N}$, be such that
$I_{h_n}\in \cl W$ and set $E = \cap_{n=1}^{\infty} h_n^{-1}(0)$.
The annihilator $\frak{M}$ of $\cl W$ in $\cl B(H_1,H_2)$ is a non-zero w*-closed bimodule.
By \cite{sht}, the support $F$ of $\frak{M}$ is not marginally null, and if an
operator $I_h$ belongs to its preannihilator then $h$ vanishes marginally
almost everywhere on $F$. It follows that
$F$ is marginally contained in $E$ and hence $E$ is not marginally null.

Now we prove the general statement.
It was shown in \cite{sht} (see the proof of \cite[Theorem~2.1]{sht})
that $\nul(\overline{\cl W}) = \supp(\cl W^{\perp})$ where $\overline{\cl W}$ is the norm
closure of $\cl W$ in $\cl C_1(H_2,H_1)$.
Hence, up to marginally null sets, $\supp(\cl W^{\perp}) \subseteq \nul(\cl W)
\subseteq \cap_{n=1}^{\infty}f_n^{-1}(0)$ and it suffices to show that $\cap_{n=1}^{\infty}f_n^{-1}(0)\subseteq
\supp(\cl W^{\perp})$. Let $\kappa$ be a rectangle disjoint from $\supp(\cl W^{\perp})$.
By the definition of support, the restriction to $\kappa$ of the functions corresponding to operators in $\cl W$
form a set dense in $\Gamma(\kappa)$. By the first part of the proof, the intersection of the null sets of the restrictions of
$f_n$, $n\in \bb{N}$, to $\kappa$ is marginally null. Hence $\kappa$ is marginally disjoint from $\cap_{n=1}^{\infty}f_n^{-1}(0)$. Since the
complement of $\supp(\cl W^{\perp})$ is $\omega$-closed, this implies the last remaining inclusion.
\end{proof}

\subsection{Schur multipliers and Peller's Theorem}

For a function $\varphi\in$ $\frak{B}(X$ $\times$ $Y)$, set
$$D(\nph) = \{k\in L^2(X\times Y) : \nph k\in L^2(X\times Y)\}.$$
We let $S_{\varphi} : D(\nph)\rightarrow L^2(X\times Y)$ be the mapping given by
$S_{\varphi}k = \varphi k$.
We identify $S_{\nph}$ with a (densely defined linear) map
on $\cl K(H_1,H_2)$ acting by the rule $S_{\varphi}(I_k) = I_{\varphi k}$.

Note that $S_{\varphi}$ depends only on the equivalence class of $\varphi$.
Taking this into account, we will sometimes say
that  a function $\varphi$ belongs to a certain class of functions, if it is equivalent to a function from this
class. When we need to make the distinction, we will write $h\in^{\sigma} \cl F$, if $h$ is equivalent
to a function from the class $\cl F$ with respect to the measure $\sigma$.

Recall that a function $\nph\in L^{\infty}(X\times Y)$ is called a
{\it Schur multiplier} if the map $S_{\varphi}$ is bounded in the operator norm,
that is, if there exists a constant $C > 0$ such that $\|S_{\nph}(I_k)\|_{\op}\leq C\|I_k\|_{\op}$,
for all $k\in L^2(X\times Y)$. If $\nph$ is a
Schur multiplier then the mapping $S_{\nph}$ has a unique weak* continuous
extension to $\cl B(H_1,H_2)$ which will still be denoted by $S_{\nph}$.

Let $\frak{S}(X,Y)$ (or $\frak{S}(X\times Y)$) be the set of all Schur multipliers; clearly, $\frak{S}(X,Y)$ is a subalgebra of $\frak{B}(X,Y)$.
The following facts follow easily from the definition of a Schur multiplier:

\begin{lemma}\label{elem}
(i) \ If $\varphi\in \frak{S}(X\times Y)$ then $\varphi|_{\alpha\times \beta} \in \frak{S}(\alpha\times\beta)$ for all
measurable subsets $\alpha\subseteq X$ and $\beta\subseteq Y$.

(ii) If $X\times Y = \cup_{p=1}^N \kappa_p$, where all $\kappa_p$ are rectangles
and $\varphi|_{\kappa_p} \in \frak{S}(\kappa_p)$ then $\varphi\in \frak{S}(X,Y)$.
\end{lemma}

Schur multipliers were first introduced by Schur in the early 20th century in
case of discrete measures $\mu$ and $\nu$. A characterisation of this particular class
of Schur multipliers was obtained by Grothendieck in \cite{Gro}. The following generalisation
for the class defined above is due to Peller \cite{peller}.

\begin{theorem}\label{peller}
Let $\nph\in L^{\infty}(X\times Y)$. The following
conditions are equivalent:

(i) \ \ $\nph$ is a Schur multiplier;

(ii) \ there exist measurable functions $a:X\to l^2$
and $b:Y\to l^2$ such that
$$\nph(x,y) = (a(x),b(y))_{l^2}, \mbox{ a.e. on } X\times Y
\mbox{ and } \sup_{x\in X}\|a(x)\|_2 \sup_{y\in Y}\|b(y)\|_2 < \infty.$$

(iii) $\nph(x,y)k(x,y)\in^{\mu\times\nu} \Gamma(X,Y)$ whenever
$k(x,y)\in\Gamma(X,Y)$.
\end{theorem}

It follows from Peller's Theorem (and can easily be seen directly)
that if the measures $\mu$ and $\nu$ are finite
then $\frak{S}(X,Y)\subseteq \Gamma(X,Y)$.

Using modern terminology, one can say that Theorem \ref{peller} identifies the algebra
$\frak{S}(X,Y)$ with the weak* Haagerup tensor product
$L^{\infty}(X,\mu)\otimes_{\whh} L^{\infty}(Y,\nu)$ (see \cite{blecher_smith}
where this tensor product was introduced).

\subsection{General facts on closable operators}

Let $\cl X$ be a Banach space. We denote by $\cl X^*$ the dual of $\cl X$.
If $\cl S\subseteq \cl X$ (resp. $\cl T\subseteq\cl X^*$), we write $\cl S^{\perp}\subseteq\cl
X^*$ (resp. $\cl T_{\perp}\subseteq\cl X$) for the annihilator (resp. preannihilator) of $\cl S$
(resp. $\cl T$).

Let $\cl Y$ be another Banach space. By an operator from $\cl X$ into $\cl Y$
we mean a linear transformation
$T : D(T)\rightarrow \cl Y$, where $D(T)$ is a
(not necessarily closed) linear subspace of $\cl X$ called the domain of $T$.
The operator $T$ is called densely defined if $D(T)$ is norm dense in $\cl X$.
We let
$$\Gr T = \{(x,Tx) : x\in D(T)\}\subseteq \cl X\oplus\cl Y$$
be the graph of $T$.
For a subset $\cl S\subseteq \cl Y\oplus\cl X$ we set $\cl S' = \{(x,y) : (y,x)\in \cl S\}$ and let $\Gr' T = (\Gr T)'$.

We recall the definition of the adjoint $T^*$ of an operator $T : D(T)\rightarrow \cl Y$.
The domain of $T^*$ is the subspace
$$D(T^*) = \{g\in \cl Y^* : \ \exists \ f\in \cl X^* \mbox{ such that } g(Tx) = f(x), \ \forall x\in D(T)\}.$$
For $g\in D(T^*)$, one lets $T^*g$ equal to $f$ where $f\in \cl X^*$ is
the functional appearing in the definition of $D(T^*)$. Note that $g\in D(T^*)$ if and only if the linear map
$x\rightarrow g(Tx)$ from $D(T)$ into $\bb{C}$ is continuous.
By the definition of the operator $T^*$, we have that $\Gr\mbox{}'(-T^*) = (\Gr T)^{\perp}$.


Recall that an operator $T$ is called {\it closable} if the norm closed hull $\overline{\Gr T}$ of
its graph is the graph of an operator. Clearly, $T$ is closable if and only if
the conditions $(x_n)_{n=1}^{\infty}\subseteq \cl X$, $y\in \cl Y$, $\|x_n\|\rightarrow 0$ and
$\|Tx_n - y\|\rightarrow 0$ imply that $y = 0$.
We call $T$ {\it w*-closable} if
the w*-closed hull $\overline{\Gr T}^{w^*}\subseteq \cl X^{**}\oplus\cl Y^{**}$
of its graph is the graph of an operator from $\cl X^{**}$ into $\cl Y^{**}$.
Here, we identify $\cl X$ and $\cl Y$ with their canonical images in their second duals.
We have that $T$ is w*-closable if and only if
the conditions $(x_{\alpha})_{\alpha}\subseteq \cl X$, $G\in \cl Y^{**}$,
w-$\lim_{\alpha} x_{\alpha} = 0$ and w*-$\lim Tx_{\alpha} = G$ imply that $G = 0$. The weak* limit is
taken with respect to the weak* topology of $\cl Y^{**}$.

In the following proposition the equivalence (iii)$\Leftrightarrow$(iv) is well-known
(see, for example, \cite[Chapter III, Section 5]{kato}); the other implications can be proved in a similar way.

\begin{proposition}\label{wstarc}
Let $T : D(T)\rightarrow \cl Y$ be a densely defined linear operator and set $\cl D = D(T^*)$. Consider the
following conditions:

(i) \ \ $T$ is w*-closable;

(ii) \ $\overline{\cl D}^{\|\cdot\|} = \cl Y^*$;

(iii) $\overline{\cl D}^{w^*} = \cl Y^*$;

(iv) \ $T$ is closable.

\noindent Then (i)$\Longleftrightarrow$(ii)$\Longrightarrow$(iii)$\Longleftrightarrow$(iv).
\end{proposition}

\section{Local Schur multipliers}

We start this section by introducing a class of functions that will play a central role
in the paper.
For brevity, let us say that a countable family of rectangles {\it covers}
$X\times Y$, or that it is {\it a covering family}, if its union is marginally equivalent to $X\times Y$.

\begin{definition}\rm
A function $\nph\in \frak{B}(X\times Y)$ will be called a {\it local Schur multiplier} if there exists a covering family
$\{\kappa_m\}_{m=1}^{\infty}$ of rectangles in $X\times Y$ such that $\nph|_{\kappa_m}\in \frak{S}(\kappa_m)$, for each $m\in \bb{N}$.
\end{definition}

The set of all local Schur multipliers on $X\times Y$ will be denoted by $\frak{S}_{\loc}(X,Y)$.

\begin{proposition}\label{loc-sp}
The set $\frak{S}_{\loc}(X,Y)$ is a subalgebra of the algebra
$C_{\omega}(X,Y)$ of all $\omega$-continuous functions.
\end{proposition}
\begin{proof} Let $\varphi\in \frak{S}_{\loc}(X,Y)$.
By the $\sigma$-finiteness of the measure spaces $(X,\mu)$ and $(Y,\nu)$, there exists
a covering family $\{\kappa_m\}_{m=1}^{\infty}$ such that
$\nph|_{\kappa_m}\in \frak{S}(\kappa_m)$,
and $(\mu\times\nu) (\kappa_m) < \infty$ for each $m\in \bb{N}$.
It follows that $\nph|_{\kappa_m} \in \Gamma(\kappa_m)$.
By \cite[Theorem 6.5]{eks}, 
$\nph|_{\kappa_m}$ is $\omega$-continuous on $\kappa_m$, $m\in \bb{N}$.
Now for an open subset $G\subseteq \mathbb{C}$, we have that
$\nph^{-1}(G) = \cup_{m=1}^{\infty} (\kappa_m\cap\nph^{-1}(G))$ is
$\omega$-open, and hence $\varphi$ is $\omega$-continuous.

It is easy to see that for two functions $\varphi, \psi $ in $\frak{S}_{\loc}(X,Y)$, one can find
a {\it common} covering family $\{\kappa_m\}_{m=1}^{\infty}$ of
rectangles on  which both $\nph$ and $\psi$ are Schur multipliers.
Since $\frak{S}(\kappa_m)$ is an algebra,
$\varphi + \psi$ and $\varphi\psi$ belong to $\frak{S}_{\loc}(X,Y)$. \end{proof}

Let $\cl V(X,Y)$ be the space of all functions $\nph\in \frak{B}(X\times Y)$
for which there exist families $\{a_i\}_{i=1}^{\infty}\subseteq \frak{B}(X)$ and
$\{b_i\}_{i=1}^{\infty}\subseteq \frak{B}(Y)$ with the properties
$$\sum_{i=1}^{\infty} |a_i(x)|^2 < \infty, \ \ \ \sum_{i=1}^{\infty} |b_i(y)|^2 <
\infty$$ for almost all $x\in X$ and $y\in Y$, and
$$\nph(x,y) = \sum_{i=1}^{\infty} a_i(x)b_i(y), \ \ \mbox{almost everywhere on } X\times Y.$$
Using a coordinate free language we may say that
$\nph\in \cl V(X,Y)$ if and only if there exists a separable Hilbert
space $H$ and measurable functions $a : X\rightarrow H$ and $b : Y\rightarrow H$ such that $\nph(x,y) =
(a(x),b(y))$ for almost all $(x,y)\in X\times Y$.

We note that $\Gamma(X,Y)\subseteq\cl V(X,Y)$ and $\frak{S}(X,Y)\subseteq \cl V(X,Y)$. Indeed,
these function spaces correspond to the cases where the functions
$\|a(\cdot)\|, \|b(\cdot)\|$ are, respectively, square integrable and essentially bounded.

\begin{lemma}\label{N-reg}
If $\nph\in \cl V(X,Y)$ then there exist families $\{X_i\}_{i=1}^{\infty}$ and
$\{Y_j\}_{j=1}^{\infty}$ of pairwise disjoint subsets of $X$ and $Y$, respectively,
such that $\nph|_{X_i\times Y_j} \in \frak{S}(X_i,Y_j)$, for all $i,j\in \bb{N}$.
We may moreover assume that $\mu(X_i) < \infty$ and $\nu(Y_j) < \infty$, for all $i,j\in \bb{N}$.
\end{lemma}
\begin{proof} Let $a : X\rightarrow \ell^2$ and $b : Y\rightarrow \ell^2$ be measurable functions
such that $\nph(x,y) = (a(x),b(y))$, for almost all $(x,y)$.
For $i,j\in \bb{N}$, set $X_i = \{x\in X: i-1 \le \|a(x)\|_2 < i\}$ and
$Y_j = \{y\in Y: j-1\le \|b(y)\|_2 < j\}$.
Then $\nph|_{X_i\times Y_j} \in \frak{S}(X_i,Y_j)$ by Theorem \ref{peller}. Partitioning
$X_i$ and $Y_j$ into subsets of finite measure, we obtain the required decompositions.
\end{proof}

Lemma \ref{N-reg} shows, in particular, that $\cl V(X,Y)\subseteq \frak{S}_{\loc}(X,Y)$.
Our aim in this section is to show that, in fact, $\cl V(X,Y) = \frak{S}_{\loc}(X,Y)$.

\begin{lemma}\label{reg}
Let $\{\kappa_m\}_{m=1}^{\infty}$ be a covering sequence of $\omega$-open sets. Then
there exist families $\{X_i\}_{i=1}^{\infty}$ and
$\{Y_j\}_{j=1}^{\infty}$ of pairwise disjoint measurable subsets of $X$ and $Y$, respectively,
such that

(i) \ $\cup_{i=1}^{\infty} X_i$ and $\cup_{j=1}^{\infty} Y_j$ have full measure, and

(ii) each rectangle $X_i\times Y_j$ is contained in a finite union of sets from $\{\kappa_m\}_{m=1}^{\infty}$.
\end{lemma}
\begin{proof}
Let us say that a subset $E\subseteq X\times Y$ is {\it mild}
if it is contained in a finite union of sets from the family $\{\kappa_m\}_{m=1}^{\infty}$. It
suffices to show that there are increasing sequences $\{A_n\}_{n=1}^{\infty}$ and $\{B_n\}_{n=1}^{\infty}$
of measurable subsets of $X$ and $Y$, respectively, such that $\cup_{n=1}^{\infty} A_n$ and
$\cup_{n=1}^{\infty}B_n$ have full measure and all rectangles of the form $A_n\times B_n$ are mild.
Indeed, the statement would then follow by setting $X_i = A_i\setminus
A_{i-1}$, $Y_j = B_j\setminus B_{j-1}$.

Since the measure spaces $(X,\mu)$, $(Y,\nu)$ are standard we may assume that $X$ and $Y$ are
equipped with $\sigma$-compact topologies with respect to which
$\mu$ and $\nu$ are regular Borel measures.
By considering increasing sequences $\{U_n\}_{n=1}^{\infty}$ and $\{V_n\}_{n=1}^{\infty}$
of compact subsets of $X$ and $Y$, respectively, we reduce the problem
to the case where $X$ and $Y$ are compact and $\mu$ and $\nu$ are finite.

We may clearly assume that each $\kappa_m$ is a rectangle.
By \cite[Lemma~3.4]{eks}, for any $\e > 0$ there exists $X_\e\subseteq X$, $Y_\e\subseteq Y$ such
that $\mu(X\setminus X_\e)<\e$, $\nu(Y\setminus Y_\e)<\e$ and $X_\e\times Y_\e$ is contained in a
finite union of rectangles $\kappa_m$. Let
$\e_n=2^{-n}$, $L_n=\cap_{k=n}^{\infty} X_{\e_k}$ and $M_n=\cap_{k=n}^{\infty}Y_{\e_k}$. Then each $L_n\times
M_n$ is mild since it is contained in $X_{\e_n}\times Y_{\e_n}$.
Furthermore, $L_n\subseteq L_{n+1}$, $M_n\subseteq M_{n+1}$,
$$\mu(X\setminus L_n)\leq\sum_{k=n}^{\infty}\mu(X\setminus X_{\e_k})<\e2^{2-n}
\text{ and } \nu(Y\setminus M_n)\leq\sum_{k=n}^{\infty}\mu(Y\setminus Y_{\e_k})<\e2^{2-n}.$$ Thus,
$\mu(X\setminus(\cup_{n=1}^\infty L_n))=0$ and $\nu(Y\setminus(\cup_{n=1}^\infty M_n))=0$
and the proof is complete.
\end{proof}

The following result may be viewed as an analogue of Lemma \ref{elem} for local multipliers.

\begin{lemma}\label{vxy}
Let $\{X_i\}_{i=1}^{\infty}$ and $\{Y_j\}_{j=1}^{\infty}$ be families of
pairwise disjoint subsets of $X$ and $Y$, respectively, such that
$X=\cup_{i=1}^{\infty} X_i$ and $Y = \cup_{j=1}^{\infty} Y_j$.
Assume that $\varphi\in \frak{B}(X\times Y)$
is such that $\varphi|_{X_i\times Y_j}\in \frak{S}(X_i,Y_j)$ for all $i,j\in \bb{N}$.
Then $\varphi\in \cl V(X,Y)$.
\end{lemma}
\begin{proof}
Let $\nph_{i,j}(x,y) = \varphi|_{X_i\times Y_j}$.
By our assumption, $\varphi_{i,j}\in \frak{S}(X_i,Y_j)$ and hence, by Theorem \ref{peller}, there exist
measurable functions $a_{i,j} : X\rightarrow \ell^2$ and $b_{i,j} : Y\rightarrow \ell^2$
such that
$\nph(x,y)=(a_{i,j}(x),b_{i,j}(y))$ for almost all $(x,y)\in X_i\times Y_j$ and
$$\alpha_{i,j}\stackrel{def}{=}\sup_{x\in X_i}\|a_{ij}(x)\|_2  \sup_{y\in Y_j}\|b_{ij}(y)\|_2 < \infty.$$
We may clearly assume that $\sup_{x\in X_i}\|a_{ij}(x)\|_2=
\sup_{y\in Y_j}\|b_{ij}(y)\|_2$. Let $H=\oplus_{i,j} H_{i,j}$, where
$H_{i,j} = \ell^2$ for all $i,j\in \bb{N}$. Considering $a_{i,j}(x)$
as a vector in $H_{i,j}$, we define a function $a : X\rightarrow H$
in the following way: if $x\in X_k$ then set $a(x) = \oplus_{i,j}
\xi_{i,j}(x)$, where $\xi_{k,j}(x) =
a_{k,j}(x)/j\sqrt{\alpha_{k,j}}$ and $\xi_{i,j}(x) = 0$ for $i\neq
k$. Similarly, we define $b : Y\rightarrow H$ by letting, for $y\in
Y_l$, $b(y) = \oplus_{i,j} \eta_{i,j}(y)$, where $\eta_{i,l}(y) =
b_{i,l}(y)/i\sqrt{\alpha_{i,l}}$ and $\eta_{i,j}(y) = 0$ if $j\neq
l$. Then for each $i$ and $x\in X_i$ we have
$$\|a(x)\|_H^2=\sum_{j=1}^{\infty} \frac{\|a_{i,j}(x)\|^2_2}{j^2\alpha_{i,j}}\leq C,$$
where $C = \sum_{j=1}^{\infty} \frac{1}{j^2}$. Similarly, we see that $\|b(y)\|_H^2\leq C$. Moreover, for
$x\in X_i$ and $y\in Y_j$, we have that
$$(a(x),b(y))_H=\frac{(a_{i,j}(x),b_{i,j}(y))_{H_{i,j}}}{ij\alpha_{i,j}}$$ and therefore
$$\nph(x,y)\chi_{X_i\times Y_j} = ij\alpha_{i,j}(a(x),b(y))_H, \text{ for almost all }(x,y)\in X_i\times Y_j.$$

The next step is to see that there exist families
$\{p_i\}_{i=1}^{\infty}$ and $\{q_j\}_{j=1}^{\infty}$ of vectors in
$\ell^2$ such that $\alpha_{i,j}=(p_i,q_j)_{\ell^2}$, $i,j\in
\bb{N}$. We note first that $|\alpha_{i,j}|\leq c_ic_j$, where
$c_i=\text{max}\{1,|\alpha_{k,l}|:k,l\leq i\}$ and
$\displaystyle\frac{\alpha_{i,j}}{jc_ic_j}=(s_i,r_j)_{\ell^2}$,
where $r_j=e_j$, $\displaystyle
s_i=\sum_j\frac{\alpha_{i,j}}{jc_ic_j}e_j$ and
$\{e_j\}_{j=1}^{\infty}$ is the standard basis of $\ell^2$. Observe
that since
$\displaystyle\sum_{j}\frac{|\alpha_{i,j}|^2}{j^2c_i^2c_j^2}\leq
\sum_j\frac{1}{j^2}<\infty$, we have that $s_i\in \ell^2$. Setting
$p_i=c_is_i$ and $q_j=jc_jr_j$, we obtain
$\alpha_{i,j}=(p_i,q_j)_{\ell_2}$. Now let $p(x)=ip_i$ if $x\in X_i$
and $q(y)=jq_j$ if $y\in Y_j$. Then
$$\nph(x,y)=(p(x),q(y))_{\ell^2} (a(x),b(y))_H=(p(x)\otimes
a(x),q(y)\otimes b(y))_{\ell^2\otimes H}$$ for almost all $(x,y)\in X\times Y$.
\end{proof}

The following  theorem is the main result of the present section.

\begin{theorem}\label{local}
Let $\nph\in \frak{B}(X\times Y)$. The following are equivalent:

(i) \ $\nph$ is a local Schur multiplier;

(ii) $\nph\in\cl V(X,Y)$.
\end{theorem}
\proof (i)$\Rightarrow$(ii) Let $\{\kappa_m\}_{m=1}^{\infty}$
be a covering family of rectangles from the definition of a local Schur
multiplier. By Lemma \ref{reg},
there exist families $\{X_i\}_{i=1}^{\infty}$ and $\{Y_j\}_{j=1}^{\infty}$
of pairwise disjoint measurable sets of $X$ and $Y$, respectively,
whose unions have full measure and
each rectangle $X_i\times Y_j$ is contained in a finite union of sets of the form $\kappa_m$.

Since $\varphi|_{\kappa_m}\in \frak{S}(\kappa_m)$ for all $m\in \bb{N}$,
it follows from Lemma \ref{elem} that $\varphi|_{X_i\times Y_j} \in \frak{S}(X_i\times Y_j)$.
An application of Lemma \ref{vxy} implies (ii).

(ii)$\Rightarrow$(i) follows from Lemma \ref{N-reg}.
\endproof

Let $\mathcal{E}$ be the class of all rectangles $\alpha\times \beta$ such that $\nph|_{\alpha\times\beta}$ is a
Schur multiplier. Let $\kappa_{\varphi}$ be the complement of the $\omega$-union of $\mathcal{E}$. Then
$\kappa_{\varphi}$ is the smallest $\omega$-closed set with the property that $\varphi$ is a local Schur multiplier on
each rectangle disjoint from it.
We call $\kappa_{\varphi}$ {\it the set of LM-singularity of $\varphi$} (LM is for "local multiplier"). It may be considered as a measure of how far
$\varphi$ is from being a Schur multiplier.
In particular, we say that $\varphi$ is {\it extremely non-Schur multiplier} if $\kappa_{\varphi} = X\times Y$.
In Section \ref{s_toeplitz}
we will give an example of an $\omega$-continuous function which is extremely non-Schur multiplier.

\section{w*-closable multipliers}\label{s_wstar}

We now introduce two classes of functions
which, along with local Schur multipliers introduced in the previous section,
are the main objects of study in the paper.
We recall that $(X,\mu)$ and $(Y,\nu)$ are fixed standard measure spaces,
$H_1 = L^2(X,\mu)$ and $H_2 = L^2(Y,\nu)$.

\begin{definition}\rm
A function $\nph\in \frak{B}(X\times Y)$ is called a {\it w*-closable}
(resp. {\it closable}) {\it multiplier} if the map $S_{\varphi}$ is
w*-closable (resp. closable), when viewed as a densely defined linear
operator on $\cl K(H_1, H_2)$.
\end{definition}

For the sake of brevity, we will sometimes call a function w*-closable (resp. closable)
if it is a w*-closable (resp. closable) multiplier. We recall that we denote by $\frak{S}_{w^*}(X,Y)$
(resp. $\frak{S}_{{\rm cl}}(X,Y)$) the set of all w*-closable (resp. closable) multipliers. 

The operator $S_{\varphi}^*$ acting on $\cl C_1(H_2,H_1)$ can be easily described. Recall that the map
$k\to I_k$ establishes an identification of $\Gamma(X,Y)$ with $\cl C_1(H_2,H_1)$
and that for $f\in \frak{B}(X\times Y)$, we write $f\in^{\mu\times\nu}\Gamma(X,Y)$ if $f$
is $\mu\times\nu$-equivalent to a function in $\Gamma(X,Y)$.

\begin{lemma}\label{adjoint}
(i) We have that
$$D(S_{\varphi}^*) = \{I_h : h\in \Gamma(X,Y) \mbox{ and } \varphi h\in^{\mu\times\nu} \Gamma(X,Y)\}.$$
In particular, $D(S_{\varphi}^*)$ is a bimodule.

(ii) For every $I_h\in D(S_{\nph}^*)$, we have $S_{\varphi}^*(I_h) = I_{\varphi h}$.
\end{lemma}
\begin{proof}
(i) For every $k\in L^2(X\times Y)$ and $h\in \Gamma(X,Y)$, we have
$\langle I_k,I_h\rangle = \int kh d(\mu\times \nu)$. It follows that
a function $h\in \frak{B}(X\times Y)$ is equivalent to a function in
$\Gamma(X,Y)$ if and only if there exists $C > 0$ such that
$$\left|\int kh d(\mu\times \nu)\right|\le
C \|I_k\|_{\op}, \text{ for all }k\in L^2(X\times Y).$$
We now have
\begin{eqnarray*}
I_h\in
D(S_{\varphi}^*)
&\Leftrightarrow &\left|\langle S_{\varphi}(I_k), I_h \rangle \right| \le
C \|I_k\|_{\op} \text{ for all } k\in D(S_{\varphi})\\
&\Leftrightarrow& \left|\int \varphi kh d(\mu\times \nu)\right|\le
C\|I_k\|_{\op} \text{ for all }k\in D(S_{\varphi}) \\
&\Leftrightarrow& \varphi h\in^{\mu\times\nu} \Gamma(X,Y),
\end{eqnarray*}
since $D(S_{\varphi})$ is dense in $L^2(X\times Y)$.

(ii) is immediate from (i).
\end{proof}

\begin{lemma}\label{divit}
Let $\nph\in \frak{B}(X\times Y)$. The following are equivalent:

(i) \ $\nph$ is a w*-closable multiplier;

(ii)  there exists a covering family $\{\kappa_m\}_{m\in \bb{N}}$ of rectangles  and functions $s_m,t_m\in
\Gamma(X,Y)$ such that $s_m(x,y)\neq 0$ m.a.e. on $\kappa_m$ and
$$\nph(x,y) = \frac{t_m(x,y)}{s_m(x,y)}, \ \ \ \text{a.e. on } \kappa_m,\ m\in \bb{N}.$$
\end{lemma}
\begin{proof} (i)$\Rightarrow$(ii) If $\nph$ is a w*-closable multiplier then, by Proposition \ref{wstarc}, the
subspace $\cl U = D(S_{\nph}^*)$ is norm dense in $\cl C_1(H_2,H_1)$.
By Lemma \ref{sequ}, there is a sequence $\{h_n\}_{n=1}^{\infty}\subseteq\Gamma(X,Y)$
with $\{I_{h_n}\}_{n=1}^{\infty}\subseteq\cl U$ such that
$\cap_{n=1}^{\infty} h_n^{-1}(0) \simeq \emptyset$. Hence,
$\cup_{n=1}^{\infty} h_n^{-1}(\mathbb{C}\setminus \{0\}) \simeq X\times Y$.
Since all $h_n$ are $\omega$-continuous, the sets
$h_n^{-1}(\mathbb{C}\setminus \{0\})$ are
$\omega$-open whence we may assume that they are countable unions of
rectangles. We conclude that $X\times Y$ is marginally equivalent to a countable union of rectangles
$\kappa_m$, $m\in \bb{N}$, such that for every $m\in \bb{N}$ there exists a function $s_m\in \cl U$
with $s_m(x,y)\neq 0$ on $\kappa_m$.

By Lemma \ref{adjoint}, there exists a function $t_m \in \Gamma(\kappa_m)$ such that $t_m\sim \nph s_m$.
Hence $\nph(x,y) = \frac{t_m(x,y)}{s_m(x,y)}$ almost everywhere on $\kappa_m$.

(ii)$\Rightarrow$(i) Since $\nph s_m \sim t_m$ and $t_m\in \Gamma(X,Y)$, Lemma
\ref{adjoint} implies that $I_{s_m}\in D(S_{\nph}^*)$. Since
$\cap_{m=1}^{\infty} s_m^{-1}(0) \simeq \emptyset$, the space
$D(S_{\nph}^*)$ is norm dense in $\cl C_1(H_2,H_1)$ by Lemma
\ref{sequ}. By Proposition \ref{wstarc}, $S_{\nph}$ is w*-closable.
\end{proof}

The following characterisation of w*-closable multipliers is the main result of this section.

\begin{theorem}\label{wstarmult}
A function $\nph\in \frak{B}(X\times Y)$ is a w*-closable multiplier if and only if
there exist functions $t,s\in \cl V(X,Y)$ such that $s(x,y)\ne 0$ marginally almost everywhere
on $X\times Y$ and $\nph(x,y) = \frac{t(x,y)}{s(x,y)}$, almost everywhere on $X\times Y$.
\end{theorem}
\begin{proof} Let $\nph\in \frak{B}(X\times Y)$ be a w*-closable multiplier. By Lemma \ref{divit},
there exists a covering family $\{\kappa_m\}_{m\in \bb{N}}$ of rectangles such that
$$\nph(x,y) = \frac{t_m(x,y)}{s_m(x,y)}, \ \ \ \text{a.e. on }
\kappa_m,$$ for some $s_m,t_m\in \Gamma(X,Y)$ with $s_n(x,y)\neq 0$
m.a.e. on $\kappa_m$.

Using Lemma \ref{N-reg} and the inclusion $\Gamma(X,Y)\subseteq \cl V(X,Y)$
we may, if necessary, partition the sets $\kappa_m$ into smaller rectangles
and assume that the functions $t_m$ and $s_m$ belong to $\frak{S}(\kappa_m)$.

By Lemma \ref{reg}, there exist families $\{X_k\}_{k=1}^{\infty}$ and
$\{Y_l\}_{l=1}^{\infty}$ of pairwise disjoint measurable subsets of $X$ and $Y$, respectively,
such that $X = \cup_{k=1}^{\infty} X_k$, $Y = \cup_{l=1}^{\infty} Y_l$ and
each $X_k\times Y_l$ is contained in a finite union of rectangles
of the form $\kappa_m$.
We show that on each rectangle $X_k\times Y_l$ the function $\varphi$ can be written in the
form $\varphi(x,y) = \frac{t_{k,l}(x,y)}{s_{k,l}(x,y)}$ where $t_{k,l},s_{k,l}\in \Gamma(X_k,Y_l)$
and $s_{k,l}(x,y)\neq 0$ marginally almost everywhere on $X_k\times Y_l$.

Indeed, $X_k\times Y_l$ is the union of a finite number of pairwise disjoint
rectangles $Z\times W$ each of which is the intersections of some rectangles of the form $\kappa_m$ and
$X_k\times Y_l$. Fix $(x,y)\in Z\times W$. On $Z\times W$
the function $\varphi$ can be written
in the form $\frac{t_0}{s_0}$, where $t_0,s_0\in \Gamma(X,Y)$.
We set
$t_{k,l}(x,y) = t_0(x,y)$ and $s_{k,l}(x,y) = s_0(x,y)$.

Now let us define functions $s$ and $t$ on $X\times Y$ by setting $t(x,y)
= t_{k,l}(x,y)$ and $s(x,y) = s_{k,l}(x,y)$ if $(x,y)\in X_k\times Y_l$.
By their definition and Theorem \ref{local},
$s$ and $t$ belong to $\cl V(X,Y)$.

Conversely, suppose that $\nph \sim t/s$ for some functions $t,s\in
\cl V(X,Y)$ with $s(x,y)\neq 0$ for every $(x,y)\in X\times Y$. By
Lemma \ref{N-reg}, $X\times Y$ can be decomposed into a countable
union of rectangles on each of which $t$ is a Schur multiplier.
Applying the same lemma to each of these rectangles, we decompose
$X\times Y$ into the union of rectangles $\kappa_m$ on each of which
both $t$ and $s$ are Schur multipliers. By the $\sigma$-finiteness of the measure spaces,
we may moreover assume that
$(\mu\times\nu)(\kappa_m) < \infty$ for each $m\in \bb{N}$. It
follows that $s|_{\kappa_m}$ and $t|_{\kappa_m}$ are equivalent to
functions from $\Gamma(\kappa_m)$. An application of Lemma
\ref{divit} now implies that $\nph$ is a w*-closable multiplier.
\end{proof}

\begin{corollary}\label{loc-w}
The set $\frak{S}_{w^*}(X,Y)$ of all w*-closable multipliers is a subalgebra
of $\frak{B}(X\times Y)$ which contains $\frak{S}_{\loc}(X,Y)$.
Moreover, every w*-closable multiplier $\nph$ is equivalent to an
$\omega$-continuous function.
\end{corollary}
\begin{proof}
The fact that the collection of all w*-closable multipliers is an algebra
follows from Theorem \ref{wstarmult} and Proposition \ref{loc-sp}.
Theorems \ref{local} and \ref{wstarmult} imply that every local multiplier is w*-closable.

Let $\nph \in \frak{B}(X\times Y)$ be a w*-closable multiplier. By Theorem \ref{divit}, $\nph =
\frac{t}{s}$ almost everywhere on $X\times Y$, where $s,t \in \cl V(X,Y)$. By Theorem \ref{local},
$t$ and $s$ are local multipliers hence they are   $\omega$-continuous by Proposition \ref{loc-sp}.

It is easy to see  that if $f$ is an $\omega$-continuous function and $g: \mathbb{C} \to \mathbb{C}$ is
continuous on an open set containing $f(X\times Y)$ then $g\circ f$ is $\omega$-continuous. Hence
$\frac{1}{s}$ is $\omega$-continuous, and since $\omega$-continuous functions form an algebra,
$\frac{t}{s}$ is $\omega$-continuous.
\end{proof}

Let $\kappa_{\nph}^{w^*}\subseteq X\times Y$ be the complement
of the $\omega$-union of the family of all rectangles
$\alpha\times\beta$ such that $\varphi|_{\alpha\times\beta}\in \frak{S}_{w^*}(\alpha,\beta)$.
The next proposition will be useful for us in the subsequent sections.

\begin{proposition}\label{p_eqlam}
Let $\nph\in \frak{B}(X\times Y)$. Then $\kappa^{w^*}_{\nph} = \nul D(S_{\nph}^*)$.
\end{proposition}
\begin{proof}
It follows from Lemma \ref{sequ} that $\varphi$ is a w*-closable
multiplier if and only if $\nul D(S_{\nph}^*) = \emptyset$.
Applying this to an arbitrary rectangle $\alpha\times \beta \subseteq X\times Y$
together with the observation that $\nul D(S_{\nph|_{\alpha\times\beta}}^*) =
(\alpha\times\beta)\cap \nul D(S_{\nph}^*)$,
we obtain that $\alpha\times\beta$ has a marginally null intersection with $\nul D(S_{\nph}^*)$
if and only if $\varphi|_{\alpha\times \beta}$ is a w*-closable multiplier. This implies our statement.
\end{proof}

It follows from Corollary \ref{loc-w} that $\kappa_{\nph}^{w^*}\subseteq \kappa_{\nph}$.
It is natural to call the functions $\nph\in \frak{B}(X\times Y)$ for which $\kappa_{\nph}^{w^*} \simeq X\times Y$
{\it extremely non-$w^*$-closable multipliers}. We have that every extremely non-w*-closable multiplier is an
extremely non-Schur multiplier.

\section{Closable multipliers}\label{s_clos}

In this section we study the class $\frak{S}_{{\rm cl}}(X,Y)$ of closable
multipliers. Let $\nph\in \frak{B}(X\times Y)$.
Recall that the transformation $S_{\varphi}$ is
defined on the linear manifold $D(\nph) = \{h\in L^2(X\times Y) : \nph h \in L^2(X\times Y)\}$
by letting $S_{\nph}h = \nph h$ and,
after identifying $L^2(X\times Y)$ with $\cl C_2(H_1,H_2)$, is
considered as a densely defined operator on the space ${\cl
K}(H_1,H_2)$ of compact operators from $H_1$ into $H_2$. The dual
space of ${\cl K}(H_1,H_2)$ is the space $\cl{C}_1(H_2,H_1)$ of
nuclear operators; we identify it with $\Gamma(X,Y)$, and, by Lemma
\ref{adjoint}, the domain  of the adjoint operator is $D_*(\varphi)
\stackrel{def}{=} D(S_{\nph}^*) = \{h\in \Gamma(X,Y): \varphi h \in^{\mu\times\nu} \Gamma(X,Y)\}$. It
follows from Proposition \ref{wstarc} that $\varphi\in
\frak{S}_{{\rm cl}}(X,Y)$ if and only if $D_*(\varphi)$ is weak*
dense in $\Gamma(X,Y)$. Equivalently, $\varphi\in
\frak{S}_{{\rm cl}}(X,Y)$ if and only if
$D_*(\varphi)_{\bot} = 0$, where $D_*(\varphi)_{\bot}$ is the set of
all compact operators $K$ such that $\langle K,h\rangle =
0$ for all $h\in D_*(\varphi)$.
Note that $D_*(\varphi)$ is a sub-bimodule of the bimodule $\Gamma(X,Y)$ over the algebras
$L^{\infty}(X,\mu)$ and $L^{\infty}(Y,\nu)$.

Let $D\subseteq \Gamma(X,Y)$ be any bimodule. Then, for all measurable sets
$\alpha\subseteq X$, $\beta\subseteq Y$ and all $h\in D$, the function
$\chi_{\alpha}(x)\chi_{\beta}(y)h(x,y)$ belongs to $D$. One can choose
the sets $\alpha$ and $\beta$ in such a way that this function is a
Schur multiplier. Indeed, if
$h(x,y) = (a(x),b(y))$ for some square integrable
Hilbert space valued functions $a$ and $b$, then
it suffices to set $\alpha = \{x: \|a(x)\| \le N\}$ and
$\beta = \{y: \|b(y)\|\le N\}$, for some $N > 0$.
Letting $N$ tend to infinity, we moreover see
that $D\cap \frak{S}(X,Y)$ is norm dense in $D$.

We will need the following proposition.

\begin{proposition}\label{p_tomi}
Let $D_1,D_2\subseteq \Gamma(X,Y)$ be weak* dense bimodules, invariant under $\frak{S}(X,Y)$.
Then $D_1\cap D_2$ is weak* dense.
\end{proposition}
\begin{proof}
We identify the predual of $\Gamma(X,Y)$ with $\cl K(H_1,H_2)$. 
Let $K\in (D_1\cap D_2)_{\perp}$ and $\theta_i\in D_i\cap \frak{S}(X,Y)$, $i = 1,2$.
By the invariance of $D_1$ and $D_2$ under $\frak{S}(X,Y)$, we have that $\theta_1\theta_2\in D_1\cap D_2$.
Thus, $\langle K, \theta_1\theta_2\rangle = 0$ and therefore
$\langle S_{\theta_1}(K),\theta_2\rangle = 0$ for all $\theta_2\in D_2\cap \frak{S}(X,Y)$.
Since $D_2\cap \frak{S}(X,Y)$ is dense in $D_2$ and $S_{\theta_1}(K)$ is a compact operator,
we have that $S_{\theta_1}(K) = 0$. Thus, $\langle K,\theta_1\rangle = 0$ for all
$\theta_1\in D_1\cap \frak{S}(X,Y)$ and hence $K = 0$.
\end{proof}

\begin{theorem}\label{alg-clos}
$\frak{S}_{{\rm cl}}(X,Y)$ is a subalgebra of $\frak{B}({X\times Y})$.
\end{theorem}
\begin{proof}
Let $\nph_1$ and $\nph_2$ be closable multipliers.
By Theorem \ref{peller} and Lemma \ref{adjoint} (i), the bimodules $D_*(\nph_1)$ and $D_*(\nph_2)$ are invariant under
$\frak{S}(X,Y)$; moreover, $D_*(\nph_1)\cap D_*(\nph_2)\subseteq D_*(\nph_1+\nph_2)$.
Propositions \ref{p_tomi} and \ref{wstarc}
imply that $\nph_1+\nph_2$ is closable.

To verify that $\frak{S}_{{\rm cl}}(X,Y)$ is closed under products, it
suffices now to show that if $\varphi$ is closable then $\varphi^2$
is closable. Let $D = D_*(\varphi) = \{h\in \Gamma(X,Y): \varphi h
\in^{\mu\times\nu} \Gamma(X,Y)\}$ and $D_0 = \{h\in \frak{S}(X,Y)\cap \Gamma(X,Y):
\varphi h \in^{\mu\times\nu} \frak{S}(X,Y)\cap \Gamma(X,Y)\}$. Then $D_0$ is dense
in $D$ and hence in $\Gamma(X,Y)$.

The product of a Schur multiplier and a closable multiplier is closable (indeed, if $w\in \frak{S}(X,Y)$,
then $D_*(\varphi)\subseteq D_*(w\varphi)$ whence $D_*(w\varphi)$ is dense).
It follows that if $h\in D_0$ then $\psi \stackrel{def}{=} \varphi^2h = \varphi(\varphi h)$ is closable.

Fix $h\in D_0$ and let
$k\in D_*(\psi)$. Then $hk\in D_*(\varphi^2)$. Hence, if
$K\bot D_*(\varphi^2)$ then $0 = (K,hk) = \langle S_h(K),k\rangle$. Since $D_*(\psi)$ is dense,
we have that $\langle K,h\rangle = S_h(K) = 0$.
Since $D_0$ is dense, $K = 0$.
Thus $D_*(\varphi^2)$ is dense and $\varphi^2$ is closable.
\end{proof}

Following the analogy with harmonic analysis initiated in \cite{arveson}  
and later pursued in \cite{froelich}, let us call an
$\omega$-closed set $E\subseteq X\times Y$ an {\it operator $M$-set}
(respectively, {\it operator $M_1$-set}) if $E$  supports a non-zero
compact operator (resp. ${\mathfrak M}_{\min}(E)$ contains a
non-zero compact operator). Clearly, every operator $M_1$-set is an
operator $M$-set. We shall show in Section \ref{s_toeplitz} that there exist
operator $M$-sets which are not operator $M_1$-sets.
We will shortly see that the property of being or not being an operator $M$- (resp. $M_1$-) set
is important for deciding whether a given function is a closable multiplier.
We hence include a consequence of Proposition \ref{p_tomi} concerning sets which are not
operator $M$- or $M_1$-sets.

\begin{proposition}\label{p_unmm1}
Let $E_1, E_2 \subseteq X\times Y$ be $\omega$-closed sets. Suppose that
$E_1$ and $E_2$ are not operator $M$-sets (resp. not operator $M_1$-sets). Then $E_1\cup E_2$ 
is not an operator $M$-set (resp. not an operator $M_1$-set).
\end{proposition}
\begin{proof}
Suppose that $E_1$ and $E_2$ are not operator $M_1$-sets. Setting $D_i = {\mathfrak M}_{\min}(E_i)_{\perp}$,
we have that $D_i$ is a weak* dense sub-bimodule of $\Gamma(X,Y)$, $i = 1,2$.
Note that, by \cite{sht},
$D_i = \{\psi\in \Gamma(X,Y) : \psi\chi_{E_i} = 0 \mbox{ m.a.e.}\}$, $i = 1,2$.
It follows that $D_i$ is invariant under $\frak{S}(X,Y)$, $i = 1,2$, and that
$$D_1\cap D_2 = \{\psi\in \Gamma(X,Y) : \psi\chi_{E_1\cup E_2} = 0 \mbox{ m.a.e.}\}.$$
By \cite{sht} again, $(D_1\cap D_2)^{\perp} = {\mathfrak M}_{\min}(E_1\cup E_2)$.
By Proposition \ref{p_tomi}, $(D_1\cap D_2)^{\perp} \cap \cl K(H_1,H_2) = \{0\}$ and hence
$E_1\cup E_2$ is not an operator $M_1$-set.

Now suppose that $E_1$ and $E_2$ are not operator $M$-sets.
Let
$$D_i = \{\psi\in \Gamma(X,Y) : \psi \mbox{ vanishes on an }\omega\mbox{-open neighbourhood of } E_i\},$$
$i = 1,2$. By \cite{sht}, $D_i^{\perp} = \mathfrak{M}_{\max}(E_i)$, $i = 1,2$.
It is clear that $D_1$ and $D_2$ are invariant under $\frak{S}(X,Y)$ and, since
$E_1$ and $E_2$ are not operator $M$-sets, $D_1$ and $D_2$ are weak* dense in $\Gamma(X,Y)$.
By Proposition \ref{p_tomi}, $D_1\cap D_2$ is weak* dense in $\Gamma(X,Y)$.
However, $D_1\cap D_2$ equals
$$\{\psi\in \Gamma(X,Y) : \psi \mbox{ vanishes on an }\omega\mbox{-open neighbourhood of } E_1\cup E_2\}$$
and hence $(D_1\cap D_2)^{\perp} = \mathfrak{M}_{\max}(E_1\cup E_2)$. Thus, $E_1\cup E_2$ is not an
operator $M$-set.
\end{proof}

In the next theorem, we relate the notions of operator $M$- and operator $M_1$-sets to
closability of multipliers.

\begin{theorem}\label{sts2}
Let $\nph\in \frak{B}(X\times Y)$.

(i) \ If $\kappa_{\nph}^{w^*}$ is not an  operator $M$-set then $\nph$ is
a closable multiplier.

(ii) \ If  $\kappa_{\nph}^{w^*}$ is an operator $M_1$-set  then $\nph$ is
not a closable multiplier.
\end{theorem}
\begin{proof}
It follows from Proposition~\ref{wstarc} that $\nph$ is closable if and only if
$D(S_{\nph}^*)^\perp$ $\cap$ $\cl K(H_1,H_2)=\{0\}$. By \cite{sht} and Proposition~\ref{p_eqlam} we have
$$\mathfrak{M}_{\min}(\kappa_{\nph}^{w^*})\subseteq D(S_{\nph}^*)^\perp\subseteq \mathfrak{M}_{\max}(\kappa_{\nph}^{w^*})$$
which clearly implies the  statement.
\end{proof}

\begin{corollary}\label{closM}
(i) \ If  $E$ is not an operator $M$-set
and if, for every marginally disjoint from $E$ rectangle $\alpha\times\beta$,
the restriction $\nph|_{\alpha\times\beta}$ is
a w*-closable multiplier, then $\nph$ is a closable multiplier.

(ii) \ If $(\mu\times \nu)(\kappa_{\nph}^{w^*})\neq 0$ then $\nph$ is not a closable multiplier.
\end{corollary}
\begin{proof}
(i) By the definition of $\kappa_{\nph}^{w^*}$,
we have that $\kappa_{\nph}^{w^*}\subseteq E$, whence $\kappa_{\nph}^{w^*}$ is not an operator $M$-set.
The claim now follows from Theorem \ref{sts2} (i).

(ii) Note that any set $E$ of non-zero measure is an
operator $M_1$-set, because it supports a non-trivial Hilbert-Schmidt
operator, and all such operators belong to
${\mathfrak M}_{\min}(E)$ \cite{arveson}. So it suffices to apply
Theorem \ref{sts2} (ii).
\end{proof}

\begin{remark}\rm
{\bf (i) } If the set $\kappa_{\nph}^{w^*}$ is synthetic then $\nph$ is a closable multiplier
if and only if $\frak{M}_{\max}(\kappa_{\nph}^{w^*})$ does not contain a non-zero compact operator.

\smallskip

\noindent {\bf (ii) } Since $\kappa_{\nph}^{w^*} \subseteq \kappa_{\nph}$, we obtain that for the closability of $\varphi$ it suffices to show that $\kappa_{\nph}$ does not support non-zero compact operators.

\smallskip

\noindent {\bf (iii)} In Section \ref{s_toeplitz} we shall construct a non-closable multiplier
$\nph$ such that
$\kappa_{\nph}^{w^*}$ is an operator $M$-set but not an operator $M_1$-set.
\end{remark}

\smallskip

\begin{example}\label{dE}\rm
Let $E\subseteq X\times Y$ be $\omega$-closed and let $\partial E$ be its $\omega$-boundary
(that is, $\partial E = E\setminus E_0$, where $E_0$ is the largest, up to
marginal equivalence, $\omega$-open set contained in $E$
\cite{eks}).
If $\nph = \chi_E$ then for each rectangle $\alpha\times\beta$ marginally contained either in $E$ or in $E^c$,
we have  that $\nph|_{\alpha\times\beta}$ is  a Schur multiplier and hence $\kappa_\nph^{w^*}$ is marginally contained in $\partial E$.
If $\partial E$ is not an operator $M$-set
then, by Theorem~\ref{sts2}, $\chi_E$ is a closable multiplier.
\end{example}

We now present our first example of a non-closable multiplier,
using a result on spectral (non)-synthesis.

\begin{example}\label{ex_52}
\emph{Let $U$ be the bilateral shift acting on the space $\ell^2({\mathbb Z})$, that is,
the operator given by $Ue_n=e_{n+1}$, $n\in \bb{N}$, where $\{e_n\}_{n\in\bb{Z}}$ is the
standard basis of $\ell^2(\bb{Z})$.
Fix $p > 2$. By \cite[Proposition~9.9]{sht2}, there exist sequences $\{a_n\}_{n\in{\mathbb
Z}}$, $\{b_n\}_{n\in{\mathbb Z}}\in \ell^2(\bb{Z})$ with $|a_n|=|b_n|$,
and an operator $X\in \cl C_p(\ell^2(\bb{Z}))$ such that
$$\sum_{n\in {\mathbb Z}} (a_nU^n)X(b_nU^{-n})=0 \ \mbox{ and } \ \sum_{n\in {\mathbb Z}} (a_nU^n)^*X(b_nU^{-n})^*\ne 0.$$
Let $W : \ell^2({\mathbb Z})\to L^2({\mathbb T})$ be the inverse Fourier transform.
Then $WUW^*$ is the operator of multiplication by $e^{it}$ and $T=WXW^*$ is an operator in $\cl C_p(L^2(\bb{T}))$
satisfying
$$\sum _{n\in {\mathbb Z}}M_{f_n}TM_{g_n}=0 \text{ and } \sum _{n\in {\mathbb Z}}M_{\overline{f}_n}TM_{\overline{g}_n}\ne 0,$$
where $M_{f_n}$ and $M_{g_n}$ are the multiplication operators by
the functions $f_n$ and $g_n$ given by $f_n(t)=a_ne^{int}$ and  $g_n(t)=b_ne^{-int}$,
respectively.
Set $d_n=a_nb_n$ and note that $\{d_n\}\in \ell^1({\mathbb Z})$.
Let $\psi(t,s) = \sum_{n\in {\mathbb Z}}f_n(t)g_n(s)=\sum_{n\in{\mathbb Z}}d_ne^{in(t-s)}$. As
$\sum_{n\in{\mathbb Z}}|f_n(t)|^2=\sum_{n\in{\mathbb Z}}|g_n(s)|^2=\sum_{n\in{\mathbb Z}}|d_n|<\infty$
for all $s,t\in \bb{T}$, Theorem \ref{peller} shows that the function $\psi$ is a Schur multiplier on $\bb{T}\times\bb{T}$
(equipped with the product Lebesgue measure).
Let $$\nph(t,s)=
\begin{cases}
\frac{\overline{\psi(t,s)}}{\psi(t,s)}& \text{if }\psi(t,s)\ne 0\\
0&\text{otherwise.}\end{cases}$$
We claim that $\nph$ is not closable.
To see this, assume that $\{T_n\}_{n=1}^{\infty}\subseteq \cl C_2(L^2(\bb{T}))$ is a sequence with
$T_n\rightarrow T$ in the operator norm.
Then $$S_\psi(T_n)\to S_\psi(T)=\sum _{n\in {\mathbb Z}}M_{f_n}TM_{g_n}=0.$$
However,
$$S_\nph(S_\psi(T_n))=S_{\overline{\psi}}(T_n)\to
S_{\overline{\psi}}(T)= \sum _{n\in {\mathbb Z}}M_{\overline{f}_n}TM_{\overline{g}_n}\ne 0.$$}
\end{example}

Example \ref{ex_52} will be considerably strengthened later:
in Proposition \ref{c_contnonc}, we will construct an $\omega$-continuous function
which is a non-closable multiplier.
On the other hand, the above example has the advantage
that it exhibits a multiplier which is not closable in $\cl C_p$, for each $p > 2$.

Let $[0,1]$ be the unit interval equipped with the Lebesgue measure,
let $\Delta = \{(x,y)\in [0,1]\times [0,1] : x < y\}$ and $\varphi =\chi_{\Delta}$ be the characteristic
function of $\Delta$. The multiplier $S_{\varphi}$ is usually called {\it the
transformer of triangular truncation} (see for example \cite{Goh-Krein}).
The following result extends the well-known fact that $S_{\varphi}$ is not a Schur multiplier.

\begin{proposition}\label{clnotwcl}
The transformer of triangular truncation is closable but  not w*-closable.
\end{proposition}
\begin{proof} We first show that $\varphi$ is closable. Let $\Lambda = \{(x,x) : x\in [0,1]\}$
be the diagonal of the unit square.
The set $\Lambda$ only supports operators of multiplication by functions in $L^{\infty}(0,1)$; in
particular, it is not an operator $M$-set.
Since the function $\varphi$ is constant on each rectangle marginally disjoint from $\Lambda$,
the claim follows from Corollary \ref{closM} (i).

To show that $\varphi$ is not w*-closable, it suffices, by Corollary \ref{loc-w},
to show that $\varphi$ is not
equivalent to an $\omega$-continuous function. Assume, towards a contradiction,
that there exists an $\omega$-continuous function $f$ such that
$f = \varphi$  almost everywhere. By Lemma~\ref{wopen}, $f=0$ m.a.e.
on $\Delta'$ and $f=1$ m.a.e. on $\Delta$.

Note that if a rectangle is marginally disjoint from $\Delta$ or
$\Delta'$ then it is marginally disjoint from $\Lambda$. It follows
that the same is true for any $\omega$-open set. Since
$f^{-1}(\mathbb{C}\setminus \{1\})$ is marginally disjoint from
$\Delta$, we obtain that $f = 1$ m.a.e on $\Lambda$. Similarly $f=0$
m.a.e. on $\Lambda$. This is a contradiction because $\Lambda$ is
not marginally null.
\end{proof}

\noindent {\bf Remark } The proof of Proposition \ref{clnotwcl} implies the following more
general statement:
Let $\Delta_1$ and $\Delta_2$ be disjoint $\omega$-open sets and $\Lambda = (\Delta_1\cup\Delta_2)^c$
be such that (a) $\Lambda$ does not support a non-zero compact operator, and (b)
for every rectangle $\kappa$, $\kappa\cap \Lambda \not\simeq \emptyset$ implies that
$\kappa\cap \Delta_i \not\simeq \emptyset$, $i = 1,2$. Then $\chi_{\Delta_1}$ is closable but not
w*-closable.

\begin{example}\rm
Let $E\subseteq X\times Y$ be an $\omega$-closed set
such that $E\setminus\partial E\not\simeq\emptyset$, where $\partial E$ is the $\omega$-boundary of $E$,
and let $\nph = \chi_E$. Then $k_\nph^{w^*}=\nul D(S_\nph^*)=\partial E$.

In fact, if a rectangle $\kappa$  is such that   $\kappa\cap \partial E \not\simeq \emptyset$ then, by the
proof of Proposition \ref{clnotwcl}, 
$\nph|_\kappa$ is not $\omega$-continuous and hence not a w*-closable multiplier, giving that $\kappa$ is not marginally  disjoint from $\kappa_\nph^{w^*}$. As $\partial E$  marginally  contains $\kappa_\nph^{w^*}$ (see Example~\ref{dE}), we obtain  $\kappa_\nph^{w^*}\simeq\partial E$.
\end{example}

Proposition~\ref{clnotwcl} shows that there exist closable multipliers which are not
$\omega$-continuous. But they are continuous in the stronger
 pseudo-topology, $\tau$, introduced in Section~\ref{sec2}.

\begin{proposition}
Any closable multiplier is $\tau$-continuous.
\end{proposition}
\begin{proof}
Let  $\nph \in \frak{S}_{{\rm cl}}(X,Y)$. If $U\subseteq {\mathbb C}$ is an open set then
$$f^{-1}(U)=(f^{-1}(U)\cap\kappa_\nph^{w^*})\cup((f^{-1}(U)\cap(\kappa_\nph^{w^*})^c).$$
Since  $\kappa_\nph^{w^*}$ is $\omega$-closed,
$(\kappa_\nph^{w^*})^c$ is marginally  equivalent to a countable
union $\cup_{i=1}^{\infty}\alpha_i\times\beta_i$ of rectangles.
Moreover, for each $i$, $f|_{\alpha_i\times\beta_i}$ is
w*-closable and hence $\omega$-continuous. This implies that
$f^{-1}(U)\cap(\alpha_i\times\beta_i)$ is marginally  equivalent to
a countable  union of rectangles and hence the same is true for
$f^{-1}(U)\cap(\kappa_\nph^{w^*})^c$. It remains to note that
$(\mu\times \nu)((f^{-1}(U)\cap\kappa_\nph^{w^*})=0$ because, by
Corollary \ref{closM}, $(\mu\times\nu)(\kappa_\nph^{w^*})=0$.
\end{proof}

\begin{remark}\rm
We note that the class of $\tau$-continuous functions is strictly
larger than that of closable multipliers; see Proposition \ref{c_contnonc}.
\end{remark}

\section{Divided differences}

Let $f$ be a continuous function on a finite or infinite open
subinterval $J\subseteq\R$. The divided difference of $f$ is the
function
$$\check f(x,y)=\frac{f(x)-f(y)}{x-y}$$
defined on $J\times J\setminus\Lambda$, where
$\Lambda=\{(x,x):x\in\R\}$. Let $\mu$ be a regular measure on
$\R$ whose support contains $J$. 
In what follows we will assume that $\mu$ is non-atomic and hence $f$
is defined almost everywhere on $J\times J$.

The property of $\check f$ being a Schur multiplier is closely
related to a kind of \lq\lq operator smoothness'' of $f$. Recall
that $f$ is called {\it operator Lipschitz (OL)} on a compact subset
$K\subseteq J$ if there exists a constant $D > 0$ such that
$$\|f(A)-f(B)\|\leq D\|A-B\|$$
for all selfadjoint operators $A$, $B$ with spectra in $K$. The
smallest constant $D$ with this property will be denoted by
$|f|_{\OL}$.

Let $O(f)$ be the union of all open subintervals $I\subseteq J$ on which
$f$ is OL. It is an open subset of $J$. Its complement will be
denoted by $E(f)$.

\begin{lemma}\label{eqdiv}
Let $I$ be a compact subset of $J$. A function $\check f$ is a Schur
multiplier on $I\times I$ if and only if $f$ is OL on $I$.
\end{lemma}
\begin{proof}
If $\check f$ is a Schur multiplier then, for
$h_1(x,y)=(x-y)h(x,y)$, we have
\begin{equation}\label{schur}
\|I_{\check f h_1}\|\leq C\|I_{h_1}\|
\end{equation}
and hence
\begin{equation}\label{ol}
\|f(A)X-Xf(A)\|\leq C\|AX-XA\|,
\end{equation}
where $A$ is the operator of multiplication by $x$ on $L^2(I,\mu)$
and $X=I_h$. By \cite[Remark 2.1, Corollary 3.6]{ks1} and
\cite[Theorem~3.4]{ks2}, $f$ is OL.

Conversely, if $f$ is OL, then (\ref{ol}) holds for each $X\in \cl
B(L^2(I,\mu))$ by \cite[Corollary 3.6]{ks1}. This implies
(\ref{schur}) for $L^2$-functions of the form
$h_1(x,y)=(x-y)h(x,y)$, where $h\in L^2(I\times I,\mu\times\mu)$.
Since functions of this form are dense in $L^2(I\times
I,\mu\times\mu)$, and since the $L^2$-norm majorizes the operator
norm, inequality (\ref{schur}) holds for all $h_1\in L^2(I\times
I,\mu\times \mu)$. This means that $\check f$ is a Schur multiplier.
\end{proof}

\begin{lemma}\label{nonin}
If $I_1$, $I_2$ are compact intervals and $I_1\cap I_2=\emptyset$
then $\check f|_{I_1\times I_2}\in \frak{S}(I_1,I_2)$.
\end{lemma}
\begin{proof}
Since $f(x)-f(y)\in \frak{S}(I_1,I_2)$, it suffices to show that
$\frac{1}{x-y}|_{I_1\times I_2}\in \frak{S}(I_1,I_2)$. Without loss of
generality we may assume that $I_1=[0,a]$, $I_2=[b,c]$ with $b>a$.
We have
$$\frac{1}{x-y}=-\sum_{n=0}^\infty \frac{x^n}{y^{n+1}}, \ \ (x,y)\in I_1\times I_2.$$
Since
$\left\|\frac{x^n}{y^{n+1}}\right\|_{\frak{S}(I_1,I_2)}\leq\frac{a^n}{b^{n+1}}$,
the series converges in $\frak{S}(I_1,I_2)$ in norm.
\end{proof}

The following theorem gives a precise description of the set of
LM-singularity for a divided difference.

\begin{theorem}\label{kf} For every continuous function $f$, we have
$$\kappa_{\check f} \simeq \{(x,x) : x\in E(f)\}.$$
\end{theorem}
\begin{proof} Write $O(f)$ as the union of a sequence of disjoint open
intervals: $O(f)=\cup_{n=1}^\infty J_n$. For each $n$, $J_n\times
J_n$ is the union of rectangles $I_k\times I_k$, where $I_k$ are
compact subintervals of $J_n$. Since, by Lemma~\ref{eqdiv}, $\check
f|_{I_k\times I_k}\in \frak{S}(I_k, I_k)$, it follows that $\check
f|_{J_n\times J_n}\in \frak{S}_{\loc}(J_n,J_n)$. Thus $(J_n\times J_n)\cap
\kappa_{\check f}\simeq\emptyset$. Furthermore, $\kappa_{\check
f}\subseteq\Lambda$ by Lemma~\ref{nonin}. It follows that, up to a
marginally null set, we have
$$\kappa_{\check f}\subseteq \Lambda\setminus(\cup_{n=1}^\infty J_n\times J_n)= \{(x,x):x\in
E(f)\}.$$

To prove the converse inclusion, it suffices to show by the
regularity of $\mu$ that if $I_1$ and $I_2$ are compact subsets of $J$ such that
$\check f|_{I_1\times I_2}\in \frak{S}(I_1,I_2)$ then $E(f)\cap I_1\cap
I_2=\emptyset$; indeed, we would then have $(I_1\times I_2)\cap\{(x,x):x\in
E(f)\}=\{(x,x):x\in E(f)\cap I_1\cap I_2\}\simeq\emptyset$.

Let $I = I_1\cap I_2$. By Lemma \ref{elem} (i), $\check f|_{I\times I}\in \frak{S}(I,I)$, 
and Lemma \ref{eqdiv} implies that $f$ is OL on $I$; therefore
$I\subseteq O(f)$ and hence $I\cap E(f)=\emptyset$.
\end{proof}

\begin{corollary}\label{div-loc}
$\check f$ is a local Schur multiplier if and only if $\mu(E(f))=0$.
\end{corollary}

It is known \cite{ks} that the class of all continuous Schur
multipliers on $X\times Y$, where $X,Y$ are compact Hausdorff spaces, very weakly
depends on the choice of Borel measures on $X$ and $Y$: it depends
only on the support of a measure. The above corollary shows that
the class of continuous local Schur multipliers essentially depends
on the choice of a measure. Indeed, a change of the measure does not
change the set $E(f)$ while the condition $\mu(E(f))=0$ need not be
preserved.

\begin{corollary}
For each $f$, the function $\check f$ is a closable Schur
multiplier.
\end{corollary}
\begin{proof} By Theorem~\ref{kf}, $\kappa_{\check f}\subseteq\{(x,x): x\in J\}$. Since the diagonal  $\{(x,x): x\in J\}$ does not support a compact operator, it follows from
Theorem~\ref{sts2} that $\check f$ is not closable.
\end{proof}

\begin{proposition}\label{wclnotlcl}
There exists a function $f:[0,1]\to{\mathbb C}$ such that $\check f$
is a Schur multiplier, $\check f\ne 0$ almost everywhere and
 $1/\check f$ is not a local Schur multiplier.
\end{proposition}

\begin{proof}
Let $M$ be a Cantor-like set of non-zero Lebesgue measure (see \cite{hs}) and
let $g$ be a continuously differentiable function which is equal to zero on $M$ and positive
otherwise. Let $f$ be its primitive function: $f'=g$. Then $f\in C^2([0,1])$
and hence  it is  operator Lipschitz \cite{DK}; by Lemma \ref{eqdiv}, $\check f$ is  a Schur 
multiplier. Since $f$ is strictly monotone, $\check f\ne 0$ almost everywhere.

The function  $1/\check{f}$ which, since $f$ is strictly monotone, is defined almost everywhere,
is not a local Schur multiplier. In fact,
assuming the converse, given $\e>0$, we can find subsets $X_\e$, $Y_\e$ of $[0,1]$ such that
$m([0,1]\setminus X_\e)<\e$,  $m([0,1]\setminus Y_\e)<\e$ and   $(\check f)^{-1}|_{X_\e\times Y_\e}$ is a
Schur multiplier. Hence $(\check f)^{-1}$ is equivalent
to an essentially bounded function. But this is impossible since
by construction $(\check f)^{-1}(x,y)$ is arbitrary large for $(x,y)$ close to $(x,x)$, $x\in M$ and since
$m(M)>0$, the set  $\{(x,y)\in X_\e\times Y_\e:|(\check f(x,y))^{-1}|>C\}$ has positive measure for all
$C>0$ and sufficiently small $\e>0$.
\end{proof}

The divided difference $\check f$ can be extended to a continuous function on $J\times J$ if
and only if $f$ is continuously differentiable. Our next aim is to
construct a continuously differentiable function $f$ such that
$\check f$ is not a Schur multiplier on each rectangle with non
marginally null intersection with $\Lambda$. For this we need an
extension of the well-known result of Farforovskaya \cite{farf} (see
also Peller \cite{peller}) which states that a continuously
differentiable function on a compact interval need not be OL.

\begin{theorem}\label{verynonOL}
There is a function in $C^1([0,1])$ which is not OL on each
subinterval of $[0,1]$.
\end{theorem}
\begin{proof}
By \cite{farf}, there exists $f\in C^1([0,1])$ which is not operator Lipschitz on $[0,1]$. Such function
$f$ can be chosen so that
\begin{equation}\label{cond}
f(0)=f'(0)=f(1)=f'(1)=0.
\end{equation}
To see this it suffices to choose a continuously differentiable non
OL function $g$ on a subinterval $I\subseteq (0,1)$ and extend it to a
continuously differentiable function $f$ on $[0,1]$ satisfying
(\ref{cond}).

Let us denote by $C^1_0$ the set of all $f\in C^{1}([0,1])$ satisfying (\ref{cond}).
Let $C_s = C^{\infty}([0,1])\cap C_0^1$. It is well-known that all functions in
$C^{\infty}$ are OL (in fact, it suffices for $f$ to have a
continuous second derivative). We claim that for each $C>0$, there
exists $g\in C_s$, such that $\|g\|_{C^1}=1$ and $|g|_{\OL} > C$.

Indeed suppose that this is not true. Since $C_s$ is dense in
$C_0^1$, each function $f\in C_0^1$ is the sum of a series
$\sum_{n=1}^{\infty}g_n$, where $g_n\in C_s$ and $\sum_{n=1}^{\infty}\|g_n\|_{C^1} <
\infty$. Our assumption gives $\sum_{n=1}^{\infty}|g_n|_{\OL} <
\infty$  which easily implies that $|f|_{\OL} < \infty$, and so $f\in OL$.
This is a contradiction because, as we know from \cite{farf}, $C^1_0$ is not contained
in the set of all Operator Lipschitz functions.

Now, by \cite{ks1}, we may state that there exist operators $A=A^*$ and
$X$ such that
$$\|g(A)X-Xg(A)\|\geq C/2\|AX-XA\|.$$
Moreover, by \cite{ks1}, $A$ and $X$ can be chosen to have finite rank.
Clearly, the interval $[0,1]$ can be replaced by an arbitrary closed
interval.

Let $\{I_n\}$  be a sequence of subintervals of $[0,1]$  such that
each subinterval $J\subseteq [0,1]$ contains at least one (and hence
infinitely many) $I_n$.

We claim that given operators of finite rank $X_1,\ldots, X_{n-1}$,
$A_1,\ldots, A_{n-1}$, where $A_i^*=A_i$, $i=1,\ldots,n-1$, and a
number $C>0$, there exist finite rank operators $A = A^*$ and $X$,
and a smooth function $g$ such that $\supp g\subseteq I_n$,
$\|g\|_{C^1([0,1])}\leq 1$, $g(A_j)=0$, $j=1,\ldots, n-1$, and
$\|[g(A),X]\|\geq C\|[A,X]\|$.

Indeed, since the spectra of all $A_j$ are finite, one can find a
subinterval $J$ of $I_n$ having empty intersection with
$\cup_{j=1}^{n-1}\sigma(A_j)$. Now by the second paragraph, there exists a
smooth function $g$ with support in $J$, such that $\|g\|_{\OL}>C$ and
$\|g\|_{C^1}=1$.  By the previous arguments this will imply the
existence of operators $A$ and $X$ with the required properties.

This allows us to construct sequences of operators $\{X_n\}$,
$\{A_n\}$, of smooth functions $\{g_n\}$ and of positive constants $\{C_n\}$ such
that
\begin{enumerate}
\item $\|g_n\|_{C^1}\leq 1$;
\item $\supp g_n\subseteq I_n$;
\item each $X_n$, $A_n$ are of finite rank and $A_n=A_n^*$;
\item $g_n(A_j)=0$ for $j<n$;
\item $\|[g_n(A_n),X_n]\|\geq C_n\|[A_n,X_n]\|$;
\item $C_n\geq 2^n (n+\sum_{j=1}^{n-1}2^{-j} |g_j|_{\OL})$.
\end{enumerate}
Let $f(t)=\sum_{j=1}^{\infty}2^{-j}g_j(t)$ so $f\in C^1([0,1])$. Let
us prove that $f$ is not OL on any subinterval $J\subseteq [0,1]$.
Assume the converse; then there exists $J\subset [0,1]$ and $C > 0$ such that
$\|[f(A),X]\|\leq C\|[A,X]\|$ for any $X$ and $A=A^*$ with
$\sigma(A)\subseteq J$. By the choice of $I_n$, given $m>0$ there
exists $n>m$ such that $I_n\subseteq J$. Therefore
$$f(A_n)=\sum _{j=1}^{\infty}2^{-j}g_j(A_n)=\sum_{j=1}^n2^{-j}g_j(A_n).$$
Since $\|[f(A_n),X_n]\|\leq C\|[A_n,X_n]\|$, we have
\begin{eqnarray*}
&&\|2^{-n}[g_n(A_n),X_n]\|\leq C\|[A_n,X_n]\|+\sum_{j=1}^{n-1}2^{-j}\|[g_j(A_n),X_n]\|\\
&&\leq (C+\sum_{j=1}^{n-1}2^{-j} |g_j|_{\OL})\|[A_n,X_n]\|.
\end{eqnarray*}
On the other hand,
 $$\|[g_n(A_n),X_n]\|\geq C_n\|[A_n,X_n]\|.$$
 Hence
 \begin{eqnarray*}
 C_n\leq 2^n (C+\sum_{j=1}^{n-1} 2^{-j} |g_j|_{\OL}).
 \end{eqnarray*}
 From condition (6) on the constant $C_n$ we get
 $2^n(n+\sum_{j=1}^{n-1}2^{-j}\|g_j\|_{\OL})\leq 2^n(C+\sum_{j=1}^{n-1}2^{-j}\|g_j\|_{\OL})$ and hence $n\leq C$ for
 every $n\in \bb{N}$,
 a contradiction.
\end{proof}

\begin{corollary}
There exists $f\in C^1([0,1])$ such that $k_{\check
f}=\{(x,x):x\in[0,1]\}$.
\end{corollary}
\begin{proof} Let $f$ be the function constructed in Theorem
\ref{verynonOL}. Then $O(f)=\emptyset$ and $E(f)=[0,1]$. The
statement now follows from Theorem~\ref{kf}.
\end{proof}

\section{Multipliers of Toeplitz type}\label{s_toeplitz}

Let $G$ be a locally compact second countable abelian group and therefore
metrisable  by \cite[8.3]{hr}. Let
$\mu=ds$ be the left invariant Haar measure on $G$. We write $L^p(G)$ for
$L^p(G,\mu)$, $p=1,2$, and denote by $C_c(G)$ the space of all continuous functions
on $G$ with compact support.
Let $\hat{G}$ be the dual group of $G$
and $A(G)$ (resp. $B(G)$) be the Fourier (resp. the Fourier-Stieltjes) algebra of $G$.
We recall that $A(G)$ is the image of $L^1(\hat{G})$ under Fourier transform.
It is well-known that $A(G)$ coincides with the family of functions
$t\mapsto \int_G f(s-t)g(s)ds = (\lambda(t)f,\bar g)$, $f$, $g\in L^2(G)$, where
$\lambda(t)f(s)=f(s-t)$. The algebra $B(G)$ is the
image under Fourier transform of the convolution algebra $M(G)$ of all
bounded Radon measures on $G$. One has $A(G)\subseteq B(G)$;
equality holds if and only if $G$ is compact. It is moreover known
that $B(G)$ coincides with  the space of all multipliers $A(G)$ (see \cite{rudin}).

For a subset $J\subseteq A(G)$, its {\it null set} is defined by
$$\nul J = \{s\in G: f(s) = 0 \text{ for all } f\in J\}.$$
Conversely, for a closed subset $E$ of $G$ we denote by
$I(E)$ (resp. $J(E)$) the space of all $f\in A(G)$ vanishing on $E$
(resp. the closed hull of the space of all $f\in A(G)$ vanishing on a neighborhood of $E$);
we have that $I(E)$ is the largest (resp. the smallest) closed ideal of $A(G)$
whose null set is equal to $E$ (see \cite{rudin}).

Let $N$ be the map sending a measurable function $f :
G\rightarrow\bb{C}$ to the function $Nf : G\times G\rightarrow
\bb{C}$ given by $Nf(s,t) = f(s - t)$. The functions of the form
$Nf$ will be called functions of Toeplitz type. It is well-known
(see, for example, \cite{bf}) that
if $f\in L^{\infty}(G)$ then $Nf$ is a Schur multiplier with respect
to Haar measure if and only if $f\in^{\mu} B(G)$. 
In this section we show that the
algebra of w*-closable multipliers of Toeplitz type coincides with
that of local Schur multipliers of Toeplitz type; if $G$ is compact
then both spaces coincide with the algebra of Schur multipliers of
Toeplitz type, that is, with $NA(G)$.

We shall start with a result relating the continuity of a function $f$ on $G$
to the $\omega$-continuity of $Nf$.

The following lemma is certainly
known but, since we were not able to find a precise reference, we include its proof for completeness.
Let $\cl O(X)$ denote the set of all open subset of a topological space $X$.

\begin{lemma}\label{mapO}
Let $X$ be a topological space and $\xi: \cl O(\mathbb{C}) \to \cl{O}(X)$ be a union preserving map such that 
$\xi(\emptyset) = \emptyset$, $\xi(\mathbb{C}) = X$
and $\xi(U\cap V) = \emptyset$ whenever $U,V\in \cl O(\bb{C})$ and $U\cap V = \emptyset$.
Then there exists a continuous function $g: X\to \mathbb{C}$ such that $\xi(U) = g^{-1}(U)$ for all $U\in \cl{O}(\mathbb{C})$.
\end{lemma}
\begin{proof}
For $t\in X$, let $O(t)$ denote the union of all $U\in \cl{O}(\mathbb{C})$ with $t\notin \xi(U)$.
Since $\xi(\bb{C}) = X$, we have that $\mathbb{C}\setminus O(t)$ is non-empty.
If it contains at least two points, say
$\lambda_1$ and $\lambda_2$, then taking disjoint open sets $U_i$ with $\lambda_i\in U_i$, $i = 1,2$,
we obtain that $t\in \xi(U_1)\cap\xi(U_2)$. This contradicts the fact that $\xi(U_1)\cap\xi(U_2)$ is empty.

We proved that $\mathbb{C}\setminus O(t) = \{\lambda\}$, for some $\lambda\in \bb{C}$.
Setting $g(t) = \lambda$, we obtain a function $g : X\rightarrow \bb{C}$. It follows from its definition
that $g^{-1}(U) = \xi(U)$, for every $U\in \cl{O}(\mathbb{C})$. Hence, $g$ is continuous.
\end{proof}

For $t\in G$, we denote by  $\Lambda_t$ the $t$-shifted diagonal:
$$\Lambda_t = \{(x,x-t): x\in G\}.$$ We say that a subset $E$ of $\Lambda_t$ is non-null in $\Lambda_t$, if $m(\{x: (x,x-t)\in E\}) > 0$.

For $W\subseteq C\times G$ set $$\pi(W) = \{t\in G: W\cap \Lambda_t \text{ is non-null in }\Lambda_t\}.$$
Clearly  $\pi(G\times G) = G$, $\pi(\emptyset)= \emptyset$ and  $\pi(W_1\cup W_2) = \pi(W_1)\cup\pi(W_2)$.

\begin{lemma}\label{steinh}  If $W$ is $\omega$-open, then $\pi(W)$ is open.
\end{lemma}
\begin{proof}
Let $s\in \pi(W)$. It follows that there exists a rectangle
$\alpha\times\beta\subseteq W$ with non-null (in $\Lambda_s$)
intersection with $\Lambda_s$. By the $\sigma$-finiteness of the measure spaces, we
may moreover assume that $\alpha$ and $\beta$ have finite measure. We now have
$m(\alpha\cap(\beta+s)) > 0$. Since the function $x\to
m(\alpha\cap(\beta+x))$ is continuous (being the convolution of the $L^2$-functions $\chi_{\alpha}$ and $\chi_{\beta}$),
$m(\alpha\cap(\beta+x))>0$ for all $x$ in a neighborhood $V$ of
$s$. Hence $V\subseteq \pi(W)$.
\end{proof}

\begin{proposition}\label{palmost}
Let $f : G\rightarrow \bb{C}$ and $\nph = Nf$. The function $\nph$ is
equivalent to an $\omega$-continuous function if and only if $f$ is equivalent to a continuous function.
Moreover, $\nph$ is $\omega$-continuous if and only if $f$ is continuous.
\end{proposition}
\begin{proof}
If $f$ is continuous then $Nf$ is continuous and hence $\omega$-continuous.
It follows easily that if $f$ is equivalent to a continuous function then
$Nf$ is equivalent to a continuous function. We hence show the converse assertions.

Let  $\psi : G\times G\rightarrow \bb{C}$ be an $\omega$-continuous function equivalent to $Nf$.
Thus, $Z \stackrel{def}{=} \{(x,y) \in G\times G : Nf(x,y) \neq \psi(x,y)\}$ is a null set.
Then $M\stackrel{def}{=} \pi(Z)$ is a null subset of $G$.
Let us say that a point $t\in G$
is {\it good} if $t\notin M$.

For $U\in \cl{O}(\mathbb{C})$, set $\xi(U) = \pi(\psi^{-1}(U))$. It follows from Lemma \ref{steinh} that $\xi$ maps $ \cl O(\mathbb{C})$ to $\cl{O}(G)$. The conditions $\xi(\mathbb{C}) = G$, $\xi(\emptyset)= \emptyset$ and  $\xi(U_1\cup U_2) = \xi(U_1)\cup\xi(U_2)$ follow from the corresponding properties of $\pi$. We have to show that $\xi$ sends disjoint sets to disjoint sets.

Note that  if $t$ is good and $t\in \xi(U)$ then $f(t)\in U$.
Indeed, $\psi(x,x-t) \in U$ for all $x$ belonging to a certain non-null set, by the definition of $\xi(U)$.
Since $t$ is good, for almost all $x\in G$, the pair $(x,x-t)$ does not belong to $Z$ hence there exists $x\in G$ such that
$\psi(x,x-t)= Nf(x,x-t) = f(t)$.

Now if $U_1\cap U_2 = \emptyset$ and $\xi(U_1)\cap \xi(U_2) \neq \emptyset$ then
$\xi(U_1)\cap \xi(U_2)$ is an open set must contain a good point $t\in G$.
But then $f(t)\in U_i$, $1=1,2$, a contradiction.

Applying Lemma \ref{mapO} we obtain a continuous function $g: G\to \mathbb{C}$ with $\xi(U) = g^{-1}(U)$,
for all $U\in \cl{O}(\mathbb{C})$. By the above argument, $g(t) = f(t)$ for all good $t$ (indeed for each $U$ containing $g(t)$ we have that $t\in g^{-1}(U) = \xi(U)$ whence $f(t)\in U$). Thus $f$ coincides almost everywhere with the continuous function $g$.

If $Nf$ is $\omega$-continuous then all points are good and $f(t) = g(t)$, for all $t\in G$.
\end{proof}

Let us say that a measurable function $f : G\rightarrow\bb{C}$
belongs (resp. almost belongs) to $A(G)$ {\it at a point} $t\in G$ if
there exist a neighborhood $U$ of $t$ and a function $g\in A(G)$
such that $f(s) = g(s)$  everywhere
(resp. almost everywhere with respect to Haar measure $\mu$) on $U$.  If $f$  belongs to $A(G)$  at each
$t\in G$ then we say that $f$ {\it locally belongs} to $A(G)$ and write  $f\in A(G)^{\loc}$. It is obvious that
$A(G)^{\loc}\subseteq C(G)$ and, using the regularity
of $A(G)$, it is easy to show that if $G$ is
compact then $A(G)^{\loc} = A(G)$. In general we
have the inclusions $A(G)\subseteq B(G)\subseteq A(G)^{\loc}$.

If $f$ almost belongs to $A(G)$ at each point $t\in G$, it is not difficult to see that
$f$ is equivalent to a function in $A(G)^{\loc}$. We recall that in this case we write $f\in^{\mu} A(G)^{\loc}$.

For a measurable function $f : G\rightarrow\bb{C}$, let
$$J_f = \{h\in A(G) : fh \in^{\mu} A(G)\}$$
and $E_f = \nul J_f$.  Clearly, $J_f$ is an ideal of $A(G)$ whence
$$J(E_f)\subseteq \overline{J_f} \subseteq I(E_f).$$

\begin{lemma}\label{eq}
Let $f : G\rightarrow \bb{C}$ be measurable. Then
\begin{equation}\label{right}
E_f=\{t\in G: f\text{ does not almost belong to }A(G)\text{ at }t\}.
\end{equation}
\end{lemma}
\begin{proof}
Let $E$ be the set in the right hand side of (\ref{right}).
If $t\in E^c$, then $f$ almost belong to $A(G)$ at $t$ and therefore
there exists a neighborhood $V$ of $t$ such that $fg$ is equivalent to a function in $A(G)$
for any $g\in A(G)$ with $V^c\subseteq \text{null }\{g\}$.
Now if $g\in A(G)$ takes the value $1$ at $t$ and $0$ on $V^c$ then $g\in J_f$ and $t\notin \text{null }g$.
Thus, $t\not\in E_f$ and hence $E_f\subseteq E$.

To see the reverse inclusion, let $t\in E$ and assume that there exists
$g\in A(G)$ such that $fg\sim h$, $h\in A(G)$ and $g(t)\ne 0$.
Then there exists a neighborhood $U$ of $t$ such that $|g(s)|> \delta>0$ for all $s\in U$.
By the regularity of $A(G)$, we can find $q\in A(G)$ such that $q(s)g(s)=1$ for all $s\in U$;
therefore $f(s)=f(s)q(s)g(s)$ for all $s\in U$. Since $fgq \sim hq$ on $U$ and $hq\in A(G)$,
the function
$f$ almost belongs  to $A(G)$ at $t$.
We obtain a contradiction giving $E\subseteq E_f$.
\end{proof}

For notational simplicity we let $\Gamma(G) = \Gamma(G,G)$. We shall
frequently use the map $P : \Gamma(G)\to A(G)$ given by
$P(f\otimes g)(t)=(\lambda(t)g,\bar f)$.
Clearly, $P$ is a surjective contraction.

For a subset $E\subseteq G$, let $E^* = \{(x,y)\in G\times G : x - y \in E\}$.

\begin{theorem}\label{th_forf}
Assume that  $G$ is a subgroup of $\R^n$ or $\mathbb T^n$, $n\in \bb{N}$.
Let $f : G\rightarrow\bb{C}$ be a measurable function, $\nph = Nf$
and $U,V\subseteq G$ be measurable sets.
The following are equivalent:

(i) \ \ $(U\times V)\cap E_f^* \simeq\emptyset$;

(ii) \ $\nph|_{U\times V}\in \frak{S}_{\loc}(U,V)$;

(iii) $\nph|_{U\times V}\in \frak{S}_{w^*}(U,V)$.
\end{theorem}
\begin{proof}
Set $E=E_f$ and note that $(U\times V)\cap E^* \simeq \emptyset$ if and only if
$(U' - V')\cap E = \emptyset$ for some $U'\subseteq U$ and $V'\subseteq V$ with
$\mu(U\setminus U') = \mu(V\setminus V') = 0$.
We claim that $f\psi\in^{\mu} A(G)$ for every $\psi\in A(G)\cap C_c(G)$ with $\supp \psi\subseteq E^c$.
Indeed, by Lemma \ref{eq}, for each $t\in E^c$
there exists a neighborhood $V_t$ of $t$ and a function $g_t\in A(G)$ such that $f \sim g_t$ on $V_t$.
Since $\text{supp} \psi$ is compact  there exists a finite set $F\subseteq G$ such that
$\text{supp} \psi\subset \cup_{t\in F}V_t$. It follows from the regularity of $A(G)$ that
there exist $h_t\in A(G)$, $t\in F$, such that $\sum_{t\in F}h_t(x)=1$ if $x\in\text{supp} \psi$ and
$h_s(x)=0$ if $x\notin  V_s$ for each $s\in F$ (see the proof of \cite[Theorem 39.21]{hr}).
Then for every $x\in G$ we have
$$f(x)\psi(x)=\sum_{t\in F}f(x)\psi(x)h_t(x)$$
and hence $f\psi\sim\sum_{t\in F}g_t\psi h_t,$
giving  $f\psi\in^{\mu} A(G)$.

(i)$\Rightarrow$(ii)
Suppose $(U\times V)\cap E^* \simeq \emptyset$ and let
$U'\subseteq U$ and $V'\subseteq V$ be measurable subsets such that
$m(U\setminus U') = m(V\setminus V') = 0$ and $(U' - V')\cap E = \emptyset$.
Since $G$ is second countable and $U' - V' \subseteq E^c$, $m(U\setminus U')=
m(V\setminus V')=0$,
we may choose increasing sequences $\{K_n\}_{n=1}^{\infty}$ and $\{L_n\}_{n=1}^{\infty}$
of compact sets such that, up to a null set, $\cup_{n=1}^{\infty} K_n = U$ and $\cup_{n=1}^{\infty} L_n = V$,
and a compact set $M_n$ such that
$K_n - L_n\subseteq M_n \subseteq E^c$. Choose, for each $n\in \bb{N}$, a function  $\psi_n\in A(G)\cap C_c(G)$ supported in $E^c$
and taking value $1$ on $M_n$.
By the previous paragraph, $f\psi_n\in^{\mu} A(G)$ and therefore $N(f\psi_n)$ is a
Schur multiplier. Thus, for each $\xi\in \Gamma(G)$, we have
$$\nph \chi_{K_n\times L_n}\xi = N(f\psi_n)\chi_{K_n\times L_n}\xi\in^{\mu\times \mu} \Gamma(G).$$
It follows that $\nph|_{K_n\times L_n}$ is a Schur multiplier and hence $\nph\in \frak{S}_{\loc}(U,V)$.

(ii)$\Rightarrow$(iii) follows from Corollary \ref{loc-w}.

(iii)$\Rightarrow$(i) We will identify $\Gamma(U,V)$ with a subset of $\Gamma(G)$ in a natural way.
Let $\psi = \nph|_{U\times V}$. By Proposition \ref{wstarc},
$D(S_{\psi}^*)$ is norm dense in $\Gamma(U,V)$. Thus, $P(D(S_{\psi}^*))$ is norm dense in $P(\Gamma(U,V))$.
By Lemma \ref{adjoint} (i), $D(S_{\psi}^*) = \{h\in \Gamma(U\times V) : \psi h\in^{\mu\times \mu} \Gamma(U,V)\}$.
Since $fP(h)=P(\nph h)\in P(\Gamma(G))\in^{\mu} A(G)$ for every $h\in D(S_{\psi}^*)$, the set
$\{P(h) : h\in \Gamma(U,V), f P(h) \in^{\mu} A(G)\}$ is dense in $P(\Gamma(U,V))$, and hence
$$P(\Gamma(U,V))\subseteq \overline{J_f}.$$
This implies that
$$E = \nul \overline{J_f} \subseteq \nul P(\Gamma(U,V)).$$

It suffices to show that
there exist subsets $U'\subseteq U$, $V'\subseteq V$ such that
$\mu(U\setminus U') = \mu(V\setminus V')=0$ and
$(U' - V')\cap \nul P(\Gamma(U,V)) = \emptyset$ since this will imply
$$U' - V' \subseteq \nul P(\Gamma(U,V))^c \subseteq E^c$$
and hence $(U\times V) \cap E^*\simeq \emptyset$.

Let $U'$, $V'$ be the sets of density points of $U$ and $V$, respectively.
Then, by the Lebesgue density theorem (see \cite{mattila}),
$\mu(U\setminus U') = \mu(V\setminus V') = 0$.
To prove the statement,
it suffices to show that $P(\chi_U\otimes\chi_V)(s) = \mu(U\cap(V+s)) > 0$
if $s\in U'-V'$. Let $s=u-v$, $u\in U'$, $v\in V'$ and assume that
$\mu(U\cap(V+s)) = \mu((U-u)\cap(V-v)) = 0$. If
$B_\e(0)$ is the closed ball with centre $0$ and radius $\e$, we then have
$$\mu((U-u)\cap B_\e(0)) + \mu((V-v)\cap B_\e(0))\leq \mu(B_\e(0)),$$
for each $\e > 0$.
Applying the Lebesgue density theorem we obtain
$$2 = 1 + 1 = \lim_{\e\to 0}\frac{\mu((U-u)\cap B_\e(0))}{\mu(B_\e(0))} + \lim_{\e\to 0}\frac{\mu((V-v)\cap B_\e(0))}{\mu(B_\e(0))}\leq 1,$$
a contradiction.

\end{proof}

\begin{remark}\rm\label{remar}
{\bf (i)} The condition that $G$ be a subgroup of $\R^n$ or $\mathbb
T^n$ is only used to prove the implication (iii)$\Rightarrow$(i),
where we appeal to the Lebesgue density theorem. The statement
remains true for more general groups (in particular, for Lie groups), for
which there is an analog of the Lebesgue theorem (see \cite{comfort-gordon, lahiri}).

\smallskip

\noindent {\bf (ii)} Taking $U = V = G$, we see that
$\nph\in \frak{S}_{w^*}(G,G)$ implies $E_f=\emptyset$ for any group $G$,
since in this case the arguments in the proof of Theorem \ref{th_forf} give
$\overline{J_f}=P(\Gamma(G))=A(G)$.
\end{remark}

We note some consequences of Theorem \ref{th_forf}.  In the next corollary,
which gives a precise description of the sets
$\kappa_{\nph}$ and  $\kappa_{\nph}^{w^*}$, we assume that $G$ satisfies the conditions of
Theorem~\ref{th_forf} (see also Remark \ref{remar} (i)).

\begin{corollary}\label{kappa}
Let $f : G\rightarrow \bb{C}$ be a measurable function and
$\nph = Nf$. Then $\kappa_{\nph} \simeq \kappa_{\nph}^{w^*} \simeq (E_f)^*$.
\end{corollary}


The following theorem shows that the set of local Schur multipliers
and that of w*-closable multipliers coincide in the class of Toeplitz
functions.

\begin{theorem}\label{toeplitz}
Let $G$ be an arbitrary second countable locally compact abelian group.
Let $f : G\rightarrow\bb{C}$ be a measurable function and $\nph=Nf$.
The following are equivalent:

(i) \ \ $f\in^{\mu} A(G)^{\loc}$;

(ii) \ $\nph$ is a local Schur multiplier;

(iii) $\nph$ is a w*-closable multiplier.

\noindent If $G$ is compact then the above statements are equivalent to

(iv) $\nph$ is a Schur multiplier.
\end{theorem}
\begin{proof}
We note that $f\in^{\mu} A(G)^{\loc}$ if and only if $E_f=\emptyset$.
The equivalence (i)$\Leftrightarrow$(ii)$\Leftrightarrow$(iii)
follows from Theorem~\ref{th_forf} and  Remark~\ref{remar}(ii).

Assume that $G$ is compact.
Then (i)$\Leftrightarrow$(iv) follows from the
equality $A(G) = A(G)^{\loc}$ and the fact that $Nf$ is a Schur
multiplier if and only if  $f\in^m A(G)$.
\end{proof}




Our next result shows that the class of $\omega$-continuous functions is strictly larger
than the class of w*-closable multipliers.

\begin{corollary} Let $G$ be a compact abelian group, $f\in C(G)\setminus A(G)$ and
$\nph = Nf$. Then $\nph$ is $\omega$-continuous but not w*-closable.
\end{corollary}
\begin{proof} Since $\nph$ is the composition of the continuous
function $g$, given by $g(s,t)=s-t$, and the continuous function $f$, it is
$\omega$-continuous. By  Theorem~\ref{toeplitz}, $\nph$ is not a
w*-closable multiplier.
\end{proof}

\begin{example}\rm
Let $\Delta=\{(x,y)\in\R:x\leq y\}$. Then $\chi_\Delta(x,y)=\chi_{(-\infty,0]}(x-y)$
is not a w*-closable multiplier, since $\chi_{(-\infty,0]}$ does not almost belong
 to $A(\R)$ at $x=0$.
\end{example}

\begin{remark}\rm
It is known that there exists a function $f\in B(\R)$ such that
$f>0$ on $\R$, but $1/f\notin B(\R)$ (\cite[\S 32]{gelfand}). Since
a function of Toeplitz type $w(s,t)=f(s-t)$, $s,t\in\R$, $f\in C(\R)$,
is a Schur multiplier if and only if $f\in B(\R)$, we obtain a positive
Schur multiplier $w$ such that $1/w$ is not a Schur multiplier.
However, $1/w$ is a local Schur multiplier, since $1/f\in
A(\R)^{\loc}$. To see this, we note that $f$ belongs locally to
$A(\R)$ and hence for each $s\in\R$ there exists a neighbourhood
$V_s$ and $g\in A(\R)$ such that $f=g$ on $V_s$. Since $g(s)\ne 0$
on $V_s$, using the regularity of $A(\R)$ one can find $h\in A(\R)$
such that $hg=1$ on $V_s$. As $h=1/f$ on $V_s$ and $s$ is arbitrary,
we have $1/f\in A(G)^{\loc}$.

Note that, by Proposition~\ref{wclnotlcl}, there exists a Schur
multiplier $w$ such that $w(s,t)\ne 0$ almost everywhere and $1/w$ is not a
local Schur multiplier.
\end{remark}


\begin{proposition}\label{c_contnonc}
There exists an $\omega$-continuous non-closable extremely-non-Schur multiplier.
\end{proposition}
\begin{proof}
Since each continuous function on $G\times G$ is
$\omega$-continuous with respect to the Haar measure,
it suffices to exhibit a
continuous function $f$ such that $Nf$ is non-closable and $\kappa_f = G\times G$. Let $G = \T$.
By \cite[Chapter II, Theorem~3.4]{katznelson}, for any set $S\subseteq \T$
of Lebesgue measure zero there exists a function $h\in C(\T)$ whose
Fourier series diverges at every point of $S$. We can choose $S$ so
that its closure is $\T$ and take the corresponding $f\in
C(\T)$. Let $\nph=Nf$.
By the
Riemann Localisation Lemma, any function which belongs to
$A(\T)$ at $x\in\T$ has a convergent Fourier series at $x$. Thus,
$\T\subseteq E_f$ and hence  $E_f = \T$.

By Corollary~\ref{kappa}, we have $\kappa_\nph\simeq\kappa_\nph^{w^*}\simeq \T^2$ and therefore $\nph$ is extremelly-non-Schur multiplier.
Moreover, applying now Proposition~\ref{p_eqlam}, we obtain $\nul D(S_\nph^*)=\T^2$ and hence $ D(S_\nph^*)=\{0\}$, showing that $\nph$ is non-closable.
\end{proof}

Now assume $G$ is compact, so that $\hat{G}$ is discrete.
Then $A(G)=\{\sum_{\chi\in \Gamma}c_\chi\chi : \sum_{\chi\in\Gamma}|c_\chi|<\infty\}$.
The space of pseudomeasures $PM(G) = A(G)^*$ can be identified with
$\ell^{\infty}(\hat{G})$ via Fourier transform:
$F\to\{\hat F(\chi)\}_{\chi\in\Gamma}$. A pseudomeasure $F\in PM(G)$ is called a {\it pseudofunction} if  $\hat F$ vanishes at infinity.
We recall  that $PM(G)$ is an $A(G)$-module with respect to the operation $fF(g)=F(fg)$, for $F\in PM(G)$, $f,g\in A(G)$, and that the support $\text{supp }F$ of a pseudomeasure $F$ is the set
$\{x\in G: fF\ne 0\text{ whenever } f(x)\ne 0, f\in A(G)\}$.

If $E$ is a closed subset of $G$ we let
$PM(E)$ (resp. $N(E)$) denote the space of all pesudomeasures supported in $E$
(resp. the weak* closed hull of the space of all measures
$\mu\in M(G)$ supported in $E$). Here, by the weak* topology we mean
the  $\sigma(PM(G), A(G))$-topology.
Clearly, $N(E)\subseteq PM(E)$. Moreover, $PM(E)$ (resp. $N(E)$) is the largest (resp. smallest) weak* closed subspace
the support of whose every element is in $E$. Moreover, $N(E) = I(E)^{\perp}$ and $PM(E) = J(E)^{\perp}$.

Recall that a closed set $E\subseteq G$ is called an {\it $M$-set} (resp. an {\it $M_1$-set})
if $PM(E)$ (resp. $N(E)$) contains a non-zero pseudofunction.
It is known that there exists an $M$-set which is not an $M_1$-set, see \cite[Section 4.6]{graham}.
In what follows  we shall give some sufficient conditions for a
function of Toeplitz type to be a closable or a non-closable multiplier, based on the
above notions.

In \cite{froelich}, Froelich studied the question of when
a given closed set supports a non-zero compact operator and a non-zero pseudo-integral compact operators.  The next result uses ideas from \cite{froelich}.
We will use the fact that the restriction of the map $N$
(given by $Nf(x,y) = f(x-y)$) to $A(G)$
takes values in the {\it Varopoulos algebra} $V(G) \stackrel{def}{=} C(G)\hat\otimes C(G)\subseteq \Gamma(G)$.

We will need a modification of the module action of $L^1(G)$ on
$V(G)$ described on p. 365 of \cite{spronk-turowska}. For $\psi\in
\Gamma(G)$ and $r\in G$, write $r\cdot \psi$ for the function given
by $r\cdot \psi(s,t) = \psi(s+r,t+r)$. If $f\in C(G)$, let $f\cdot
\psi= \int_{G} (r\cdot \psi) f(r)dr$, where the
integral is understood in the Bochner sense. Following the arguments
in \cite{spronk-turowska}, one can show that the action on $C(G)$ on
$\Gamma(G)$ extends to an action of $L^1(G)$ on $\Gamma(G)$ and that
if $\{f_{\alpha}\}_{\alpha}$ is a bounded approximate identity for
$L^1(G)$ then $f_{\alpha}\cdot \psi\rightarrow  \psi$ for every
$\psi\in \Gamma(G)$.

The following lemma establishes that  $E$ is an $M_1$-set if and only if $E^*$ is an operator $M_1$-set. This justifies the terminology introduced in Section 5.

\begin{lemma} \label{m1set}
Let $E\subseteq G$ be a closed set. The space
${\mathfrak M}_{\min}(E^*)$ contains a non-zero compact operator if and only if  $E$ is an $M_1$-set.
\end{lemma}
\begin{proof}
Let $K$ be a non-zero compact operator supported on $E^*$. Then
there exist $\gamma,\delta\in \hat{G}$ such that $c_{\delta,\gamma}
\stackrel{def}{=} (K\gamma,\delta)\ne 0$. Let $F$ be the
pseudomeasure given by
$$\hat{F}(\chi)=c_{\gamma-\delta+\chi,\chi}, \ \ \chi\in \hat{G}.$$
Since $K$ is compact, $\hat{F}$ is a pseudo-function.
For each $v=\sum_{\chi\in \hat{G}}a_{\chi}\chi \in I(E)$
(the sum being absolutely convergent), we have
\begin{eqnarray*}F(v)&=&\sum_{\chi\in \hat{G}}a_{\chi}\hat{F}(\chi) = \sum_{\chi}a_{\chi}(K\chi,\gamma - \delta + \chi)\\
& = & \langle K,\sum_{\chi\in \hat{G}} a_{\chi}\chi\overline{\chi(\gamma - \delta)}\rangle=
\langle K, \tilde{v}\rangle,
\end{eqnarray*}
where $\langle \cdot,\cdot\rangle$ is the duality between $\cl B(L^2(G))$ and $\Gamma(G)$
and $\tilde{v}$ is the function given by $\tilde{v}(s,t) = Nv(s,t) \overline{(\gamma - \delta)(t)}$.
Since $Nv$ vanishes on $E^*$ and $K\in {\mathfrak M}_{\min}(E^*)$, we have that $F(v)=0$,
showing that $F$ is a pseudofunction in $N(E)$. Thus $E$ is an $M_1$-set.

Conversely, assume that $E$ is an $M_1$-set and let $F$ be a non-zero pseudofunction in $N(E)$.
We let $K$ be the operator on $L^2(G)$ defined by
$K\chi = \hat{F}(\chi)\chi$ on the orthonormal basis $\hat{G}$ of $L^2(G)$.
Then for $v\in I(E)$, we have
$$\langle K, Nv\rangle=F(v)=0.$$

Suppose that $\psi\in \Gamma(G)$ vanishes marginally almost everywhere
on $E^*$. For $\chi\in \hat{G}$, define the
functions $\psi^\chi$ and $\tilde \psi^\chi$ by
$$\psi^\chi(s,t) = \chi\cdot \psi(s,t) \text{ and } \tilde \psi^\chi(s,t)=\chi(s)\psi^\chi(s,t).$$
We have that $\psi^\chi, \tilde \psi^\chi\in \Gamma(G)$, and $\tilde
\psi^\chi(s+r,t+r)=\tilde \psi^\chi(s,t)$ marginally almost
everywhere, for each $r\in G$ (see \cite[Theorems~3.1 and
4.6]{spronk-turowska}). Therefore, by
\cite[Proposition~4.5]{spronk-turowska}, $\tilde \psi^\chi\in
NA(G)$. Since $\tilde \psi^{\chi}$ vanishes on $E^*$, $\tilde
\psi^\chi=Nv$ for some $v$ vanishing on $E$. By the previous paragraph, $\langle K,
\tilde \psi^\chi\rangle=0$. This implies that  $\langle K,
\psi^\chi\rangle=0$. In fact, if $\chi\equiv 1$, this is trivial; if
$\chi\not\equiv 1$ we have  $\psi^\chi(s,t)=\chi(-s)\tilde
\psi^\chi(t,s)=\chi(-s) v^{\chi}(s-t)$ for some $v^\chi\in A(G)$ and
then writing $v^\chi(s)=\sum_{\tau\in\hat G} a_\tau\tau(s)$ and taking into account that 
$(\bar\chi\tau,\tau) = 0$ for all $\tau\in \hat{G}$,
we obtain
$$\langle K,\psi^\chi\rangle=\sum_{\tau\in \hat{G}} a_\tau(K\bar\chi\tau,\tau)=
\sum_{\tau\in \hat{G}} a_\tau F(\bar\chi\tau)(\bar\chi\tau,\tau) = 0.$$

Finally, we let $\{u_\alpha\}$ be a bounded approximate identity for
$L^1(G)$ chosen from $\text{span}\{\chi : \chi\in \hat G\}$. For
each $\alpha$, we have that $u_\alpha\cdot \psi\in
\text{span}\{\psi^\chi:\chi\in \hat G\}$ and hence $\langle K,
u_\alpha\cdot \psi\rangle =0$ giving $\langle K,\psi\rangle=0$.
Since $\psi$ is an arbitrary element of $\Gamma(G)$ vanishing on
$E^*$, we conclude that $K\in {\mathfrak M}_{\min}(E^*)$.
\end{proof}

\begin{proposition}\label{clmm1}
Let $f : G\rightarrow\bb{C}$ be a measurable function and $\nph =
Nf$. Then the following holds:

(i) \ \ If $E_f$ is not an $M$-set then $\nph$ is  closable.

(ii) \ If $E_f$ is an $M_1$-set then $\nph$ is not closable.
\end{proposition}
\begin{proof}
(i) By \cite[Theorem 1.2.7, Lemma 1.2.10]{froelich}, $\frak{M}_{\max}(E_f^*)$ does not
contain a non-zero compact operator. The statement now follows from
Theorem \ref{sts2} (i) and Corollary \ref{kappa}.

(ii) follows from Theorem \ref{sts2} (ii), Corollary \ref{kappa} and Lemma \ref{m1set}.
\end{proof}

\begin{remark}
We note that  if $E_f$ satisfies spectral synthesis then $\nph$ is closable
if and only if $E_f$ is an $M$-set.
\end{remark}


We say that a subset $E\subseteq G$ is $\tau_0$-open if $E$ is equivalent,
with respect to the Haar measure, to an open subset of $G$.
A function $f:G\to{\mathbb C}$ is said to be {\it $\tau_0$-continuous}
if $f^{-1}(U)$ is $\tau_0$-open for any open $U\subseteq G$.

\begin{proposition}
Let $f : G\rightarrow\bb{C}$ be a measurable function and $\nph =
Nf$. If $\nph$ is closable then $f$ is $\tau_0$-continuous.
\end{proposition}
\begin{proof}
Since $f$ almost belongs to $A(G)$ at each point $t\in E_f^c$,
it is equivalent to a continuous function $h$ on $E_f^c$. In fact,
for each $t\in E_f^c$ there exists a neighborhood  $U_t$ and $h_t\in A(G)$ such that $f=h_t$ almost everywhere on $U_t$. Since $G$ is second countable there exists a countable number of neighborhoods $U_{t_i}$ such that $E_f^c=\cup_{i=1}^{\infty} U_{t_i}$. Now set
$h(t)=h_{t_i}(t)$  for $t\in U_{t_i}$. Clearly, $h$ is continuous on $E_f^c$ and $f=h$ a.e. on $E_f^c$.
Thus, given an open subset $U\subseteq G$, the set $f^{-1}(U)\cap E_f^c$ is equivalent to an open set in $G$. Since
$$f^{-1}(U)=(f^{-1}(U)\cap E_f^c)\cup(f^{-1}(U)\cap E_f),$$
it is enough to show that $\mu(E_f)=0$.

Assume, by way of contradiction, that $\mu(E_f) > 0$.
Since $G$ is a compact group,
$\chi_{E_f}\in L_1(G)$, $\hat{\chi}_{E_f}$ 
vanishes at infinity and  hence $\chi_{E_f}$ is a non-zero pseudo-function supported in $E_f$. But since, by Proposition~\ref{clmm1}, $E_f$ is
not an $M_1$-set, we arrive at a contradiction.
\end{proof}

We will finish this section by constructing
an example  of a non-closable  multiplier $\nph$ for which $\kappa_\nph^{w^*}$
is an operator $M$-set but not an operator $M_1$-set.
The existence of a closable multiplier for which $\kappa_\nph^{w^*}$
is an operator $M$-set but not an operator $M_1$-set remains an open problem.



\begin{example}\rm
Let $E\subseteq\T$ be s an $M$-set which is not an $M_1$-set. Then
$\partial E = E$. We have that $\mathfrak{M}_{\min}(E^*)$ does not contain a
non-zero compact operator, while $\mathfrak{M}_{\max}(E^*)$ contains
such an operator, say $K$.

As $\mathfrak{M}_{\max}(E^*)\ne \mathfrak{M}_{\min}(E^*)$, we can find $\Psi\in \Gamma(\bb{T})$ which vanishes on $E^*$
such that $\langle K,\Psi\rangle\ne 0$.

Let $\displaystyle\Psi_1=\sum_n\frac{1}{2^n}\chi_{\alpha_n}\otimes\chi_{\beta_n}$, where $\{\alpha_n\times\beta_n\}$ is a  disjoint family of rectangles such that $(E^*)^c\simeq\cup_n\alpha_n\times\beta_n$.
Then $\Psi_1$ is the limit of elements of $\Gamma(\bb{T})$ which vanish on an $\omega$-open subset
of $\bb{T}\times\bb{T}$ containing $E^*$, and hence $\langle K,\Psi_1\rangle=0$ \cite{sht}.
We note that, moreover, $\text{null}\Psi_1=E^*$.

Now let
$$\nph(x,y)=\left\{\begin{array}{cl}\frac{\Psi(x,y)}{\Psi_1(x,y)}& (x,y)\in (E^*)^c,\\
0,& (x,y)\in E^*.\end{array}\right.$$

As $\Psi_1\in \Gamma(\T)$, one can find measurable subsets $K_n$, $n\in \bb{N}$, with
$K_n\subseteq K_{n+1}$, $n\in \bb{N}$, such that
$m(K_n^c)\to_{n\rightarrow\infty} 0$ and $\Psi_1\chi_{K_n\times K_n}, \Psi\chi_{K_n\times K_n}\in \frak{S}(K_n,K_n)$, $n\in \bb{N}$.
Then there exists $N$ such that $\langle K,\Psi\chi_{K_N\times K_N}\rangle\ne 0$.
On the other hand, we have $S_{\Psi_1\chi_{K_N\times K_N}}(K) = 0$.
Let $T_n\in \cl C_2(L^2(\T))$, $T_n\to K$, $n\to \infty$.
Then $S_n\stackrel{def}{=}M_{\chi_{K_N}}T_n M_{\chi_{K_N}}\to M_{\chi_{K_N}}KM_{\chi_{K_N}}$ and
$$S_{\Psi_1}(S_n)=S_{\Psi_1\chi_{K_N\times K_N }}(S_n)\to S_{\Psi_1\chi_{K_N\times K_N }}(K)=0$$
but
$$S_\nph(S_{\Psi_1}(S_n))=S_\Psi(S_n)\to S_\Psi( M_{\chi_{K_N}}KM_{\chi_{K_N}})\ne 0.$$

Thus $\nph$ is not closable and hence $\kappa_\nph^{w^*}$ is an operator $M$-set.
As $\kappa_\nph^{w^*}\subseteq E^*$, we have
${\mathfrak M}_{\min}(\kappa_\nph^{w^*})\subseteq {\mathfrak M}_{\min}(E^*)$ and hence $\kappa_\nph^{w^*}$ is not an operator $M_1$-set.
\end{example}

\section{ Open problems}

In this section we list some open problems.
The most important question which we have left unanswered is the following:

\bigskip

\noindent
{\bf Problem 1}. Is every w*-closable multiplier a local Schur multiplier?

\bigskip

\noindent
{\bf Problem 2}. Does Theorem \ref{th_forf} hold for all locally compact abelian groups?

\bigskip

\noindent
{\bf Problem 3}. For which $f\in C(\mathbb{R})$ is the divided difference $\check{f}$ a w*-closable multiplier?

\medskip

Since all local multipliers are w*-closable, Corollary
\ref{div-loc} shows that a sufficient condition for this to happen is $\mu(E(f))=0$.

\bigskip

The last two problems are related to Problem 1.

\bigskip

\noindent
{\bf Problem 4}. Let $f$ be an Operator Lipschitz function on $[a,b]$, and let $f^{\prime}(x)\neq 0$ for all $x\in [a,b]$. Is the inverse function $f^{-1}$ Operator Lipschitz on $[f(a),f(b)]$?

\bigskip

\noindent
{\bf Problem 5}. For which continuous functions $f$ and normal operators $A\in \cl B(H)$
is the map on $\cl B(H)$ given by $AX-XA \to f(A)X-Xf(A)$ closable in weak* topology?

\bigskip

\noindent {\bf Remarks (i) }
The  map considered in Problem 5 is norm closable.
Indeed, if $AX_n-X_nA\to_{n\rightarrow \infty} 0$ and
$f(A)X_n-X_nf(A)\to_{n\rightarrow \infty} B$ then
$$[B,A] = \lim_{n\rightarrow \infty} [[f(A),X_n],A] =
\lim_{n\rightarrow \infty} [f(A),[X_n,A]] = 0,$$
that is, $B$ belongs to the commutant $\{A\}^{\prime}$ of $A$. If
$\cl E : \cl B(H) \to \{A\}^{\prime}$ is a conditional expectation,
then $B = \cl E(B) = \lim_{n\rightarrow \infty} \cl E(f(A)X_n - X_nf(A)) = \lim_{n\rightarrow \infty} (f(A)\cl E(X_n) - \cl E(X_n)f(A)) = 0$.

\smallskip

\noindent {\bf (ii)}
For the case $f(z) = \overline{z}$ the answer to Problem 5 is negative.
More precisely the \lq\lq Fuglede'' map $AX - XA \to A^*X - XA^*$ is not w*-closable,
if $\sigma(A)$ has non-empty interior and the spectral measure of $A$ is
equivalent to the Lebesgue measure on the interior $U$ of $\sigma(A)$.

To see this, we assume for simplicity that $A$ is the operator of
multiplication by $z$ on $L^2(U,dz d\overline{z})$. Let $f$ be the function on $G = \bb{R}^2$
given by $f(z) = \frac{\overline{z}}{z}$ and let
$\varphi = Nf$ be the corresponding Toeplitz multiplier on $G\times G$.
It is not difficult to
check that the set $E_f$ of all points $s\in G$ at which $f$ does not belong to
$A(G)$ is the singleton $\{0\}$; applying Corollary \ref{kappa} we get that
$\kappa_{\varphi}^{w^*} = \Lambda = \{(z,z): z\in G\}$. It follows
that the multiplier $\varphi$ is not w*-closable on $U\times U$, so there
are Hilbert Schmidt operators $I_{h_n}$ supported in $U\times U$ with $I_{h_n}\to 0$
and $I_{\varphi h_n} \to B \neq 0$ in the weak* topology. We may assume that $h_n(z_1,z_2)$
vanish on some neighborhoods of the diagonal $\Lambda$. Indeed, let $V_n$ be a
neighborhood of $\Lambda$ such that the Hilbert-Schmidt norm $\|h_n\chi_{V_n}\|_2$ is less than $1/n$.
Then $\|\varphi h_n\chi_{V_n}\|_2 < 1/n$ whence $\|\varphi h_n\chi_{V_n}\| < 1/n$ and
we may replace $h_n$ by $h_n - h_n\chi_{V_n}$.

Setting $p_n(z_1,z_2) = h_n(z_1,z_2)/(z_1-z_2)$ and $X_n = I_{p_n}$ we get that
$[A,X_n] = I_{h_n} \to 0$ and $[A^*,X_n] = I_{\varphi h_n} \to B$.


\begin{thebibliography}{99}


\bibitem{arveson}
{\sc W.B. Arveson}, {\it Operator algebras and invariant subspaces},
{\rm Ann. Math. (2) 100 (1974), 433-532}

\bibitem{BS1} \textsc{M.S. Birman and M.Z. Solomyak}, \textit{Stieltjes
double-integral operators. II}, \textrm{(Russian) Prob. Mat.Fiz. 2
(1967), 26-60}

\bibitem{BS2}\textsc{M.S. Birman and M.Z. Solomyak}, \textit{Stieltjes
double-integral operators, III (Passage to the limit under the
integral sign)}, \textrm{(Russian) Prob. Mat. Fiz. No 6 (1973),
27-53}

\bibitem{BS3} \textsc{M.S. Birman and M.Z. Solomyak}, \textit{Operator
Integration, perturbations and commutators}, \textrm{Zap. Nauchn.
Sem. Leningrad. Otdel. Mat. Inst. Steklov. (LOMI) Issled. Linein.
Oper. Teorii Funktsii. 17, 170 (1989), 34-66}

\bibitem{BS4}
\textsc{M.S. Birman and M.Z. Solomyak}, \textit{ Double operator
integrals in a Hilbert space},  \textrm{Int. Eq. Oper. Th. 47
(2003),  no. 2, 131--168}

\bibitem{blecher_smith}
{\sc D.P. Blecher and R. Smith}, {\it The dual of the Haagerup
tensor product}, \textrm{J. London Math. Soc. (2) 45 (1992),
126--144}

\bibitem{bf} \textsc{M. Bozejko and G. Fendler},
\textit{Herz-Schur multipliers and completely bounded multipliers of
the Fourier algebra of a locally compact group}, \textrm{Colloquium
Math. 63 (1992) 311-313}

\bibitem{comfort-gordon} {\sc W.W. Comfort and H. Gordon}, {\it Vitali's theorem for invariant measures},  {\rm Trans. Amer. Math. Soc.  99  1961 83--90}


\bibitem{DK} {\sc J.L. Daletskii and S.G. Krein}, \textit{Integration and
differentiation of functions of hermitian operators and applications to the
theory of perturbations}, Amer. Math. Soc. Translations (2) 47 (1965), 1-30

\bibitem{eks}
{\sc J.A. Erdos, A. Katavolos and V.S. Shulman},
{\it Rank one subspaces of Bimodules over Maximal Abelian Selfadjoint Algebras},
{\rm J. Funct. Anal. 157 No.2 (1998), 554-587}

\bibitem{farf}
{\sc Yu.B. Farforovskaya}, {\it An estimate of the norm $||f(A)-f(B)||$ for selfadjoint operators $A$ and $B$}, {\rm Zap. Nauchn. Semin. LOMI 56 (1976), 143--162 (English transl. J. Sov. Math. 14 (1980), 1133--1149)}

\bibitem{froelich} {\sc J. Froelich}, {\it Compact operators, invariant subspaces and spectral synthesis}, {\rm J. Funct. Anal. 81 (1988), 1-37}

\bibitem{gelfand} {\sc I. Gelfand, D. Raikov, G. Shilov}, {\it Commutative normed rings. Translated from the Russian, with a supplementary chapter},
    {\rm Chelsea Publishing Co., New York, 1964}

\bibitem{Goh-Krein} {\sc I.C. Gohberg and M.G. Krein}, {\it Theory and applications of Volterra operators in Hilbert space},
{\rm Translation of Mathematical Monographs, vol. 24, American Mathematical Society, 1970}


\bibitem{graham}{\sc C.Graham and O.C. McGehee} {\textit Essays in commutative harmonic analysis},
{\rm Springer-Verlag, New York-Berlin, 1979}


\bibitem{Gro} \textsc{A. Grothendieck}, \textit{Resume de la theorie metrique
des produits tensoriels topologiques}, \textrm{Boll. Soc. Mat.
Sao-Paulo 8 (1956), 1-79}

\bibitem{hr}\textsc{E. Hewitt and K.A. Ross}, \textit{Abstract harmonic analysis. Vol. I. Structure of topological groups, integration theory, group representations}, \textrm{Springer-Verlag, Berlin-New York, 1979}

\bibitem{hs} \textsc{E. Hewitt and K. Stromberg}, \textit{Real and abstract analysis. A modern treatment of the theory of functions of a real variable}, \textrm{Springer-Verlag, New York, 1965}





\bibitem{kato} {\sc T. Kato},
{\it Perturbation theory for linear operators}, {\rm Springer-Verlag, Berlin, 1995}



\bibitem{katznelson} {\sc Y. Katznelson, } {\it An introduction to harmonic analysis Third edition}, {\rm Cambridge University Press, Cambridge, 2004}

\bibitem{ks} {\sc E. Kissin and V.S. Shulman,}
{\it Operator multipliers}, {\rm Pacific J. Math. 227  (2006),  no. 1, 109--141}

\bibitem{ks1} {\sc E. Kissin and V.S. Shulman} {\it Classes of operator-smooth functions. I. Operator-Lipschitz functions.} {\rm  Proc. Edinb. Math. Soc. (2)  48  (2005),  no. 1, 151--173}

\bibitem{ks2} {\sc E. Kissin, and V.S. Shulman,} {\it On the range inclusion of normal derivations: variations on a theme by Johnson, Williams and Fong.} {\rm  Proc. London Math. Soc. (3)  83  (2001),  no. 1, 176--198}

\bibitem{lahiri} {\sc B.K. Lahiri},
{\it On translations of sets in topological groups}
{\rm J. Indian Math. Soc. (N.S.) 39 (1975), 173--180}

\bibitem{mattila} {\sc P. Mattila},{\it  Geometry of sets and measures in Euclidean spaces. Fractals and rectifiability}, {\rm Cambridge University Press, Cambridge, 1995}



\bibitem{peller} {\sc V. Peller,} {\it Hankel operators in the theory of perturbations of unitary and selfadjoint operators. (Russian)}  {\rm Funktsional. Anal. i Prilozhen.  19  (1985),  no. 2, 37--51, 96}


\bibitem{sht} {\sc V.S. Shulman and L. Turowska,} {\it Operator synthesis. I. Synthetic sets, bilattices and tensor algebras}, {\rm J. Funct. Anal.
209 (2004), 293-331}

\bibitem{sht2}{\sc V.S. Shulman and L. Turowska,}{\it Operator synthesis II: Individual synthesis and linear operator equations} {\rm  J. Reine Angew. Math.  590  (2006), 143--187}

\bibitem{spronk-turowska}{\sc N. Spronk and L. Turowska,} {\it Spectral synthesis and operator synthesis for compact groups},
{\rm J. London Math. Soc. (2) 66 (2002), 361-376}

\bibitem{rudin} {\sc W. Rudin,}  {\it Fourier analysis on groups}, {\rm John Wiley \& Sons, Inc., New York, 1990}

\end{thebibliography}
\end{document}